\documentclass[11pt,a4paper]{article}

\usepackage{setspace}
\usepackage[T1]{fontenc}
\usepackage[utf8]{inputenc}
\usepackage{authblk}
\usepackage{tikz}
\usepackage{subfigure}
\usepackage{pgfplots}
\usetikzlibrary{decorations.markings}
\usetikzlibrary{plotmarks}

\usepackage{graphics}
\usepackage{float}
\usepackage{amssymb}
\usepackage{amsfonts}
\usepackage{amsmath}
\usepackage{amsthm}
\usepackage{pdfsync}
\newtheorem{theorem}{Theorem}
\newtheorem{lemma}[theorem]{Lemma}
\newtheorem{definition}[theorem]{Definition}

\usepackage{array}
\usepackage{IEEEtrantools}

\makeatletter
\newcommand{\thickhline}{%
    \noalign {\ifnum 0=`}\fi \hrule height 1pt
    \futurelet \reserved@a \@xhline
}
\newcolumntype{"}{@{\hskip\tabcolsep\vrule width 1pt\hskip\tabcolsep}}
\makeatother

\pgfplotsset{tick label style={font=\Huge},    label style={font=\Huge}, legend style={font=\Huge}}

\usepackage{amssymb}

\usepackage{amsthm}

\usepackage[normalem]{ulem}

\newcommand{\ie}{{\it i.e.},\ }
\newcommand{\eg}{{\it e.g.},\ }

\title{Necessary conditions for the invariant measure of a random walk to be a sum of geometric terms}
\author[1]{Yanting Chen}
\author[1]{Richard J. Boucherie}
\author[1,2]{Jasper Goseling}
\affil[1]{Stochastic Operations Research, University of Twente,The Netherlands}
\affil[2]{Department of Intelligent Systems, Delft University of Technology,The Netherlands}
\affil[ ]{\small\texttt{\{y.chen,r.j.boucherie,j.goseling\}@utwente.nl}}

\date{\today}

\begin{document}

\maketitle

\begin{abstract}
We consider the invariant measure of homogeneous random walks in the quarter-plane. In particular, we consider measures that can be expressed as an infinite sum of geometric terms. We present necessary conditions for the invariant measure to be a sum of geometric terms. We demonstrate that, under a mild regularity condition, each geometric term must individually satisfy the balance equations in the interior of the state space.  We show that the geometric terms in an invariant measure must be the union of finitely many pairwise-coupled sets of infinite cardinality. We further show that for the invariant measure to be a sum of geometric terms, the random walk cannot have transitions to the north, northeast or east. Finally, we show that for an infinite weighted sum of geometric terms to be an invariant measure at least one of the weights must be negative.
\end{abstract}



\section{Introduction} \label{sec:introduction}
Random walks for which the invariant measure is a geometric product form are often used to model practical systems. The benefit of using such models is that their performance can be readily analyzed with tractable closed-form expressions. However, the class of random walks which have a product form invariant measure is rather limited. Therefore, it is of interest to find larger classes of tractable measures that can be the invariant measures for random walks in the quarter-plane. 

We study random walks in the quarter-plane that are homogeneous in the sense that transition probabilities are translation invariant. Our interest is in finite measures $m(i,j)$ that can be expressed as
\begin{equation}{\label{eq:invariantmeasure}}
 m(i,j) = \sum_{(\rho,\sigma) \in \Gamma} \alpha(\rho, \sigma)\rho^i\sigma^j, \quad \text{with} \quad |\Gamma| = \infty,
\end{equation}
\ie $m(i,j)$ can be expressed as an infinite sum of geometric terms. 

Examples of random walks with invariant measure of form~\eqref{eq:invariantmeasure} exist. One such an example is the $2 \times 2$ switch, which is introduced in~\cite{boxma1993compensation}. In~\cite{boxma1993compensation}, the joint distribution is presented explicitly for the $2 \times 2$ switch problem in the form of a sum of two alternating series of product-form geometric distributions.

Contrary to much other work, for instance~\cite{cohen1983boundary,fayolle1999random,miyazawa2011light,neuts1981matrix}, our interest is not in finding the invariant measure for specific random walks. Instead our interest is in characterizing the fundamental properties of random walks, sets $\Gamma$ and coefficients $\alpha$ that allow an invariant measure to be expressed in form~\eqref{eq:invariantmeasure}. 
In~\cite{chen2012invariant} we investigated invariant measures that are a sum of finitely many geometrics for the random walk in the quarter-plane. The results that are presented in the current paper form the natural extension of our results from~\cite{chen2012invariant} to the case of countably many terms.

In the current work we demonstrate that under a mild technical condition the following necessary conditions must hold:
\begin{enumerate}
\item The geometric measures that correspond to the elements of $\Gamma$ individually satisfy the balance equations in the interior of the state space.
\item Set $\Gamma$ is the union of finitely many sets, each of which have infinite cardinality and have a structure that we refer to as pairwise-coupled.
\item In the interior of the state space, the random walk has no transitions to the north, northeast or east.
\item At least one of the coefficients $\alpha(\rho,\sigma)$ in~\eqref{eq:invariantmeasure} is negative. 
\end{enumerate}


Adan et al.~\cite{adan1993compensationXXX,adan1991analysis,adan1993compensation,adan1990symmetric} use a compensation approach to construct an invariant measure that is an infinite sum of pairwise-coupled geometric terms. In particular, this is done for the random walk that have no transitions to the north, northeast and east. In addition the geometric terms in the measure will individually satisfy the interior balance equations and negative weights will occur in the sum. The work of Adan et al.\  demonstrates that measures that satisfy the necessary conditions obtained by us may indeed be an invariant measure of a random walk.

The balance equation in the interior of the state space induce an algebraic curve. A closely related curve arises as the kernel of the boundary value problems studied in~\cite{cohen1983boundary,fayolle1999random} and related work. Some of its basic properties were derived in~\cite{fayolle1999random}. An important part of the current work consists of studying this algebraic curve in more detail. Section~\ref{sec:propertyQ} presents new results on the geometric properties of this curve. 

A related study for the reflected Brownian motion in a wedge was performed by Dieker et al.~\cite{dieker2009reflected}. It was shown in~\cite{dieker2009reflected} that for the invariant measure of this process to be a linear combination of finitely many exponential measures, there must be an odd number of terms that have a pairwise-coupled structure. The methods that are developed in~\cite{dieker2009reflected} for the continuous state space Brownian motion, however, cannot be used for the discrete state space random walk.

The remainder of this paper is structured as follows. In Section~\ref{sec:model} we present the model. The geometric properties of the algebraic curve arising from the balance equations in the interior of the state space are studied in Section~\ref{sec:propertyQ}. The main contributions of the paper, necessary conditions for the invariant measure of a random walk to be a sum of geometric terms are given in Section~\ref{sec:CandS}. In Section~\ref{sec:switch} we provide an example of a random walk for which the invariant measure satisfies all the conditions that are given in Section~\ref{sec:CandS}. Section~\ref{sec:conclusion} provides concluding remarks.
\section{Model} \label{sec:model}
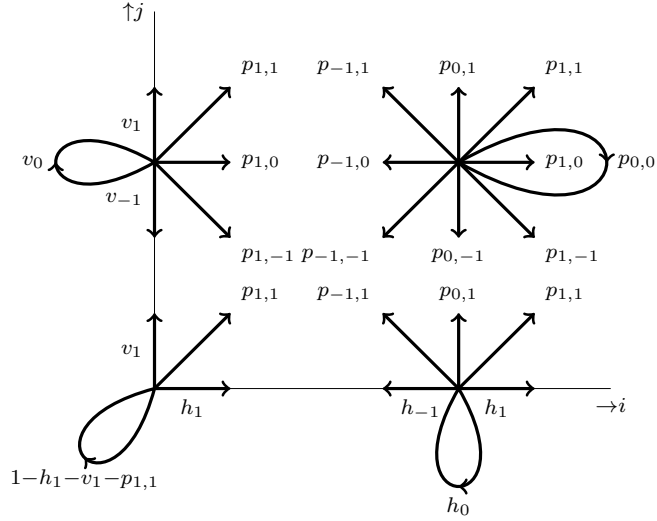
\begin{figure}
\hfill
\begin{tikzpicture}[scale=1]
\tikzstyle{axes}=[very thin] \tikzstyle{trans}=[very thick,->]
\tikzstyle{intloop}=[->,to path={
    .. controls +(30:3) and +(-30:3) .. (\tikztotarget) \tikztonodes}]
\tikzstyle{hloop}=[->,to path={
    .. controls +(-60:2) and +(-120:2) .. (\tikztotarget) \tikztonodes}]
\tikzstyle{vloop}=[->,to path={
    .. controls +(-150:2) and +(-210:2) .. (\tikztotarget) \tikztonodes}]
\tikzstyle{oloop}=[->,to path={
    .. controls +(255:2) and +(195:2) .. (\tikztotarget) \tikztonodes}]

   \draw[axes] (0,0)  -- node[at end, below] {$\scriptstyle \rightarrow i$} (6,0); 
   \draw[axes] (0,0) -- node[at end, left] {$\scriptstyle {\uparrow} {j}$} (0,5);
   \draw[trans] (0,0) to node[below] {$\scriptstyle h_{1}$} (1,0);
   \draw[trans] (0,0) to node [at end, anchor=south west] {$\scriptstyle p_{1,1}$} (1,1);
   \draw[trans] (0,0) to node[left] {$\scriptstyle v_{1}$} (0,1);
   \draw[trans] (4,0) to node[below] {$\scriptstyle h_{-1}$} (3,0);
   \draw[trans] (4,0) to node[below] {$\scriptstyle h_{1}$} (5,0);
   \draw[trans] (4,0) to node[at end, anchor = south west] {$\scriptstyle p_{1,1}$} (5,1);
   \draw[trans] (4,0) to node[at end, anchor = south]  {$\scriptstyle p_{0,1}$} (4,1);
   \draw[trans] (4,0) to node[at end, anchor = south east] {$\scriptstyle p_{-1,1}$} (3,1);
   \draw[trans] (0,3) to node[left] {$\scriptstyle v_{-1}$} (0,2);
   \draw[trans] (0,3) to node[at end, anchor = west] {$\scriptstyle p_{1,0}$} (1,3);
   \draw[trans] (0,3) to node[left] {$\scriptstyle v_{1}$} (0,4);
   \draw[trans] (0,3) to node[at end, anchor = north west] {$\scriptstyle p_{1,-1}$} (1,2);
    \draw[trans] (0,3) to node[at end, anchor = south west] {$\scriptstyle p_{1,1}$} (1,4);
   \draw[trans] (4,3) to node[at end,anchor = west] {$\scriptstyle p_{1,0}$} (5,3);
   \draw[trans] (4,3) to node[at end,anchor = south west] {$\scriptstyle p_{1,1}$} (5,4) ;
   \draw[trans] (4,3) to node[at end, anchor = south] {$\scriptstyle p_{0,1}$} (4,4);
   \draw[trans] (4,3) to node[at end, anchor = south east] {$\scriptstyle p_{-1,1}$} (3,4);
   \draw[trans] (4,3) to node[at end, anchor = east] {$\scriptstyle p_{-1,0}$} (3,3);
   \draw[trans] (4,3) to node[at end, anchor = north east] {$\scriptstyle p_{-1,-1}$} (3,2);
   \draw[trans] (4,3) to node[at end, anchor = north] {$\scriptstyle p_{0,-1}$} (4,2);
   \draw[trans] (4,3) to node[at end, anchor = north west] {$\scriptstyle p_{1,-1}$} (5,2);  
\draw[very thick,decoration={markings, mark=at position 0.5 with {\arrow{>}}},postaction={decorate}] (0,0) to[oloop] node[below] {$\scriptstyle 1 - h_1 - v_1 - p_{1,1} $} (0,0);
\draw[very thick,decoration={markings, mark=at position 0.5 with {\arrow{>}}},postaction={decorate}] (4,0) to[hloop] node[below] {$\scriptstyle h_0$} (4,0);
\draw[very thick,decoration={markings, mark=at position 0.5 with {\arrow{>}}},postaction={decorate}] (4,3) to[intloop] node[right] {$\scriptstyle p_{0,0}$} (4,3);
\draw[very thick,decoration={markings, mark=at position 0.5 with {\arrow{>}}},postaction={decorate}] (0,3) to[vloop] node[left] {$\scriptstyle v_0$} (0,3);
\end{tikzpicture}
\hfill{}
\caption{Random walk in the quarter-plane. \label{fig:rw}}
\end{figure}

Consider a two-dimensional random walk $R$ on the pairs $S = \{(i,j), i,j \in \mathbb{N}_{0}\}$ of non-negative integers. We refer to $\{(i,j) | i>0, j>0\}$, $\{(i,j) | i>0, j=0\}$, $\{(i,j) | i=0, j>0\}$ and $(0,0)$ as the interior, the horizontal axis, the vertical axis and the origin of the state space, respectively. The transition probability from state $(i,j)$ to state $(i+s,j+t)$ is denoted by $p_{s,t}(i,j)$. Transitions are restricted to the adjoining points (horizontally, vertically and diagonally), \ie $p_{s,t}(k,l)=0$ if $|s|>1$ or $|t|>1$. The process is homogeneous in the sense that for each pair $(i,j)$, $(k,l)$ in the interior (respectively on the horizontal axis and on the vertical axis) of the state space
\begin{equation} \label{eq:homogeneous}
p_{s,t}(i,j)=p_{s,t}(k,l)\quad\text{and}\quad p_{s,t}(i-s,j-t)=p_{s,t}(k-s,l-t),
\end{equation}
for all $-1\leq s\leq 1$ and $-1\leq t\leq 1$.  We introduce, for $i>0$, $j>0$, the notation $p_{s,t}(i,j)=p_{s,t}$, $p_{s,0}(i,0)=h_s$ and $p_{0,t}(0,j)=v_t$.
Note that the first equality of~\eqref{eq:homogeneous} implies that the transition probabilities for each part of the state space are translation invariant. The second equality ensures that also the transition probabilities entering the same part of the state space are translation invariant.
The above definitions imply that $p_{1,0}(0,0)=h_1$ and $p_{0,1}(0,0)=v_1$. The model and notations are illustrated in Figure~\ref{fig:rw}. 

We assume that the random walk has invariant measure $m$, \ie for $i, j > 0$,
\begin{align}
m(i,j) &= \sum_{s = -1}^{1} \sum_{t = -1}^{1} m(i-s,j-t)p_{s,t}, \label{eq:interior} \\
m(i,0) &= \sum_{s = -1}^{1} m(i-s,1) p_{s,-1} + \sum_{s = -1}^{1} m(i-s,0)p_{s,0}, \label{eq:horizontal} \\
m(0,j) &= \sum_{t = -1}^{1} m(1,j-t) p_{-1,t} + \sum_{t = -1}^{1} m(0,j-t)p_{0,t}. \label{eq:vertical} 
\end{align}
We will refer to the above equations as the balance equations in the interior, the horizontal axis and the vertical axis of the state space respectively. The balance at the origin is implied by the balance equation for all other states.

We are interested in measures that are linear combinations of geometric measures. We first classify the geometric measures.
\begin{definition}[{Geometric measures}]{\label{de:degenerateGT}}
The measure $m(i,j) = \rho^i \sigma^j$ is called a geometric measure. It is called horizontally degenerate if $\sigma = 0$, vertically degenerate if $\rho = 0$ and non-degenerate if $\rho > 0$ and $\sigma > 0$. We define $0^0 \equiv 1$.
\end{definition}

We represent a geometric measure $\rho^i \sigma^j$ by its pair of parameters $(\rho, \sigma)$ in $[0, \infty)^2$. Then, a set $\Gamma \subset [0, \infty)^2$ characterizes a set of geometric measures. Let $\mathbb{R}_{+}^{2} = \{(x,y)| x \geq 0, y \geq 0\}$, $R_B = \{(x,y) \in \mathbb{R}_{+}^{2} |xy = 0\}$, $U = \{(x,y)|(x,y) \in (0,1)^2\}$ and $\bar{U} = [0,1)^2$.

\begin{definition}[{Induced measure}]\label{def:induced}
The measure $m$ is called induced by $\Gamma\subset\mathbb{R}_{+}^2$ if
\begin{equation*}
m(i,j) = \sum_{(\rho, \sigma) \in \Gamma} \alpha(\rho, \sigma)\rho^i\sigma^j,
\end{equation*}
with $\alpha(\rho, \sigma) \in \mathbb{R}\backslash \{0\}$ for all $(\rho, \sigma) \in \Gamma$.
\end{definition}

\begin{figure}
\hfill
\subfigure[\label{fig:2a}]
{
%
%
%
%
\begin{tikzpicture}[scale = 0.2]

\begin{axis}[%
view={0}{90},
width=4.52083333333333in,
height=4.52083333333333in,
scale only axis,
xmin=0, xmax=1.8,
ymin=0, ymax=1.8,
xtick = {-0.5,0,0.5,1,1.5},
ytick = {-0.5,0.5,1,1.5}]

\addplot [draw=red, ultra thick] coordinates{ (0.363636363636364,1.4144017180055)(0.363599457655072,1.41414141414141)(0.358293676367411,1.37373737373737)(0.353552859618717,1.33333333333333)(0.349358154758861,1.29292929292929)(0.345695843347829,1.25252525252525)(0.342557224085144,1.21212121212121)(0.339938618106187,1.17171717171717)(0.337841502046863,1.13131313131313)(0.336272780717225,1.09090909090909)(0.335245218869528,1.05050505050505)(0.334778061114233,1.01010101010101)(0.334897881517189,0.96969696969697)(0.335639721217989,0.929292929292929)(0.337048595630851,0.888888888888889)(0.33918148557815,0.848484848484849)(0.342109973926918,0.808080808080808)(0.345923758612449,0.767676767676768)(0.350735376670035,0.727272727272727)(0.356686632351321,0.686868686868687)(0.363636363636364,0.648281072912731)(0.364069735766821,0.646464646464647)(0.375897026996241,0.606060606060606)(0.390046547588157,0.565656565656566)(0.404040404040404,0.532457127946925)(0.407926181867272,0.525252525252525)(0.434470377019749,0.484848484848485)(0.444444444444444,0.47244094488189)(0.47244094488189,0.444444444444444)(0.484848484848485,0.434470377019749)(0.525252525252525,0.407926181867272)(0.532457127946925,0.404040404040404)(0.565656565656566,0.390046547588157)(0.606060606060606,0.375897026996241)(0.646464646464647,0.364069735766821)(0.648281072912731,0.363636363636364)(0.686868686868687,0.356686632351321)(0.727272727272727,0.350735376670035)(0.767676767676768,0.345923758612449)(0.808080808080808,0.342109973926918)(0.848484848484849,0.33918148557815)(0.888888888888889,0.337048595630851)(0.929292929292929,0.335639721217989)(0.96969696969697,0.334897881517189)(1.01010101010101,0.334778061114233)(1.05050505050505,0.335245218869528)(1.09090909090909,0.336272780717225)(1.13131313131313,0.337841502046863)(1.17171717171717,0.339938618106187)(1.21212121212121,0.342557224085144)(1.25252525252525,0.345695843347829)(1.29292929292929,0.349358154758861)(1.33333333333333,0.353552859618718)(1.37373737373737,0.358293676367411)(1.41414141414141,0.363599457655072)(1.4144017180055,0.363636363636364)(1.45454545454545,0.373575250058153)(1.4949494949495,0.385077253152639)(1.53535353535354,0.398291070708535)(1.55118959054279,0.404040404040404)(1.57575757575758,0.427361969563184)(1.59103880063035,0.444444444444444)(1.58599755101201,0.484848484848485)(1.57575757575758,0.498430462672308)(1.55728453777577,0.525252525252525)(1.53535353535354,0.545805957620431)(1.5158327969253,0.565656565656566)(1.4949494949495,0.582456241811882)(1.46770300443741,0.606060606060606)(1.45454545454545,0.615927537936713)(1.41648436554555,0.646464646464647)(1.41414141414141,0.64817532256267)(1.37373737373737,0.679081653560027)(1.36412783641646,0.686868686868687)(1.33333333333333,0.710211591536339)(1.31204966190001,0.727272727272727)(1.29292929292929,0.741942607491584)(1.26105534769676,0.767676767676768)(1.25252525252525,0.774377885750849)(1.21212121212121,0.807598844810139)(1.21155401326555,0.808080808080808)(1.17171717171717,0.841399360965841)(1.16358095833415,0.848484848484849)(1.13131313131313,0.876430954405704)(1.11743486973948,0.888888888888889)(1.09090909090909,0.912768293900369)(1.07311749366256,0.929292929292929)(1.05050505050505,0.950504218137117)(1.03057927276445,0.96969696969697)(1.01010101010101,0.989740278713788)(0.989740278713788,1.01010101010101)(0.96969696969697,1.03057927276445)(0.950504218137117,1.05050505050505)(0.929292929292929,1.07311749366256)(0.912768293900369,1.09090909090909)(0.888888888888889,1.11743486973948)(0.876430954405704,1.13131313131313)(0.848484848484849,1.16358095833415)(0.841399360965841,1.17171717171717)(0.808080808080808,1.21155401326555)(0.807598844810139,1.21212121212121)(0.774377885750849,1.25252525252525)(0.767676767676768,1.26105534769676)(0.741942607491584,1.29292929292929)(0.727272727272727,1.31204966190001)(0.710211591536339,1.33333333333333)(0.686868686868687,1.36412783641646)(0.679081653560027,1.37373737373737)(0.64817532256267,1.41414141414141)(0.646464646464647,1.41648436554555)(0.615927537936713,1.45454545454545)(0.606060606060606,1.46770300443741)(0.582456241811882,1.4949494949495)(0.565656565656566,1.5158327969253)(0.545805957620431,1.53535353535354)(0.525252525252525,1.55728453777577)(0.498430462672308,1.57575757575758)(0.484848484848485,1.58599755101201)(0.444444444444444,1.59103880063035)(0.427361969563185,1.57575757575758)(0.404040404040404,1.55118959054279)(0.398291070708535,1.53535353535354)(0.385077253152639,1.4949494949495)(0.373575250058153,1.45454545454545)(0.363636363636364,1.4144017180055)(NaN,NaN)};
\addplot [
color=black,
solid,
line width = 1.0pt
]
coordinates{
 (0,1)(1,1) 
};

\addplot [
color=black,
solid,
line width = 1.0pt
]
coordinates{
 (1,0)(1,1) 
};

\end{axis}
\end{tikzpicture}
} 
\subfigure[\label{fig:2b}]
{
%
%
%
%
\begin{tikzpicture}[scale = 0.2]

\begin{axis}[%
view={0}{90},
width=4.52083333333333in,
height=4.52083333333333in,
scale only axis,
xmin=0, xmax=1.4,
ymin=0, ymax=1.4,
xtick = {-0.5,0,0.5,1,1.4},
ytick = {-0.5,0.5,1,1.4}]

\addplot [draw=red, ultra thick] coordinates{ (0.444444444444444,0.855072463768116)(0.4427595386728,0.848484848484849)(0.433506785021937,0.808080808080808)(0.425294122263819,0.767676767676768)(0.418464019350937,0.727272727272727)(0.413527834034268,0.686868686868687)(0.411284562799714,0.646464646464647)(0.413056181651223,0.606060606060606)(0.421199256508875,0.565656565656566)(0.440363228242016,0.525252525252525)(0.444444444444444,0.520430107526882)(0.484848484848485,0.502659574468085)(0.525252525252525,0.503770739064857)(0.565656565656566,0.51387326584177)(0.596857619295741,0.525252525252525)(0.606060606060606,0.529194187582563)(0.646464646464647,0.550306658878505)(0.674585775928315,0.565656565656566)(0.686868686868687,0.573511967064406)(0.727272727272727,0.60031847133758)(0.735808047776655,0.606060606060606)(0.767676767676768,0.631121778387336)(0.787544666332545,0.646464646464647)(0.808080808080808,0.665022421524664)(0.832729455703381,0.686868686868687)(0.848484848484849,0.703258145363409)(0.872043317184383,0.727272727272727)(0.888888888888889,0.747530864197531)(0.905994778322112,0.767676767676768)(0.929292929292929,0.800298661765551)(0.934970443852366,0.808080808080808)(0.958526801448963,0.848484848484849)(0.96969696969697,0.872869318181818)(0.977144955787214,0.888888888888889)(0.99021223901417,0.929292929292929)(0.997881606342452,0.96969696969697)(0.999463517547583,1.01010101010101)(0.994145764871894,1.05050505050505)(0.980964925573209,1.09090909090909)(0.96969696969697,1.11226851851852)(0.957295434623459,1.13131313131313)(0.929292929292929,1.16213674830968)(0.917751097303632,1.17171717171717)(0.888888888888889,1.19086021505376)(0.848484848484849,1.20862649128174)(0.833792470156105,1.21212121212121)(0.808080808080808,1.21738505747126)(0.767676767676768,1.21865430296786)(0.727272727272727,1.21284829721362)(0.725119533896965,1.21212121212121)(0.686868686868687,1.19908998657894)(0.646464646464647,1.17744593881857)(0.639081121252997,1.17171717171717)(0.606060606060606,1.14600113442995)(0.59218447097235,1.13131313131313)(0.565656565656566,1.1031746031746)(0.556985406466246,1.09090909090909)(0.528574445417293,1.05050505050505)(0.525252525252525,1.04569687738005)(0.506784273547258,1.01010101010101)(0.485864447179728,0.96969696969697)(0.484848484848485,0.967647058823529)(0.470608349396229,0.929292929292929)(0.456063907044299,0.888888888888889)(0.444444444444444,0.855072463768116)(NaN,NaN)};
\addplot [
color=black,
solid,
line width = 1.0pt
]
coordinates{
 (0,1)(1,1) 
};

\addplot [
color=black,
solid,
line width = 1.0pt
]
coordinates{
 (1,0)(1,1) 
};

\end{axis}
\end{tikzpicture}
}
\subfigure[\label{fig:2c}] 
{
%
%
%
%
\begin{tikzpicture}[scale = 0.2]

\begin{axis}[%
view={0}{90},
width=4.52083333333333in,
height=4.52083333333333in,
scale only axis,
xmin=0, xmax=1.4,
ymin=0, ymax=1.4,
xtick = {-0.5,0,0.5,1,1.4},
ytick = {-0.5,0.5,1,1.4}]

\addplot [draw=red, ultra thick] coordinates{ (0.202020202020202,0.353682214701604)(0.198384884176416,0.323232323232323)(0.197077895125914,0.282828282828283)(0.202020202020202,0.246508855441227)(0.203210994526056,0.242424242424242)(0.242424242424242,0.203210994526056)(0.246508855441227,0.202020202020202)(0.282828282828283,0.197077895125914)(0.323232323232323,0.198384884176416)(0.353682214701604,0.202020202020202)(0.363636363636364,0.203895429726532)(0.404040404040404,0.214494599816697)(0.444444444444444,0.225785196573688)(0.484848484848485,0.236997817356168)(0.504404788147808,0.242424242424242)(0.525252525252525,0.250700917059662)(0.565656565656566,0.266180848525009)(0.606060606060606,0.28040464748359)(0.613128657834267,0.282828282828283)(0.646464646464647,0.298405096645817)(0.686868686868687,0.315717377758959)(0.705478529499738,0.323232323232323)(0.727272727272727,0.334879782200734)(0.767676767676768,0.355000115222371)(0.786170830398225,0.363636363636364)(0.808080808080808,0.376913130462884)(0.848484848484849,0.399657695376722)(0.856667825816326,0.404040404040404)(0.888888888888889,0.426152914840109)(0.917591296514448,0.444444444444444)(0.929292929292929,0.453964973403925)(0.969690941847002,0.484848484848485)(0.96969696969697,0.484854384740536)(1.01010101010101,0.522600209342491)(1.01309284479204,0.525252525252525)(1.04839442169988,0.565656565656566)(1.05050505050505,0.568865091941115)(1.07584234204985,0.606060606060606)(1.09090909090909,0.636320777503951)(1.09610259029538,0.646464646464647)(1.10924344886468,0.686868686868687)(1.11585119923513,0.727272727272727)(1.11615563320778,0.767676767676768)(1.11040146355617,0.808080808080808)(1.09884567781367,0.848484848484849)(1.09090909090909,0.867567246365329)(1.08131783270928,0.888888888888889)(1.05778644774189,0.929292929292929)(1.05050505050505,0.939687761385969)(1.02741480688696,0.96969696969697)(1.01010101010101,0.989122644265788)(0.989122644265788,1.01010101010101)(0.96969696969697,1.02741480688696)(0.939687761385969,1.05050505050505)(0.929292929292929,1.05778644774189)(0.888888888888889,1.08131783270928)(0.867567246365329,1.09090909090909)(0.848484848484849,1.09884567781367)(0.808080808080808,1.11040146355617)(0.767676767676768,1.11615563320778)(0.727272727272727,1.11585119923513)(0.686868686868687,1.10924344886468)(0.646464646464647,1.09610259029538)(0.636320777503952,1.09090909090909)(0.606060606060606,1.07584234204985)(0.568865091941115,1.05050505050505)(0.565656565656566,1.04839442169988)(0.525252525252525,1.01309284479204)(0.522600209342491,1.01010101010101)(0.484854384740536,0.96969696969697)(0.484848484848485,0.969690941847002)(0.453964973403925,0.929292929292929)(0.444444444444444,0.917591296514448)(0.426152914840109,0.888888888888889)(0.404040404040404,0.856667825816326)(0.399657695376722,0.848484848484848)(0.376913130462884,0.808080808080808)(0.363636363636364,0.786170830398225)(0.355000115222371,0.767676767676768)(0.334879782200734,0.727272727272727)(0.323232323232323,0.705478529499738)(0.315717377758959,0.686868686868687)(0.298405096645817,0.646464646464647)(0.282828282828283,0.613128657834268)(0.28040464748359,0.606060606060606)(0.266180848525009,0.565656565656566)(0.250700917059662,0.525252525252525)(0.242424242424242,0.504404788147808)(0.236997817356168,0.484848484848485)(0.225785196573688,0.444444444444444)(0.214494599816697,0.404040404040404)(0.203895429726532,0.363636363636364)(0.202020202020202,0.353682214701604)(NaN,NaN)};
\addplot [
color=black,
solid,
line width = 1.0pt
]
coordinates{
 (0,1)(1,1) 
};

\addplot [
color=black,
solid,
line width = 1.0pt
]
coordinates{
 (1,0)(1,1) 
};

\end{axis}
\end{tikzpicture}
}
\subfigure[\label{fig:2d}]
{
\input{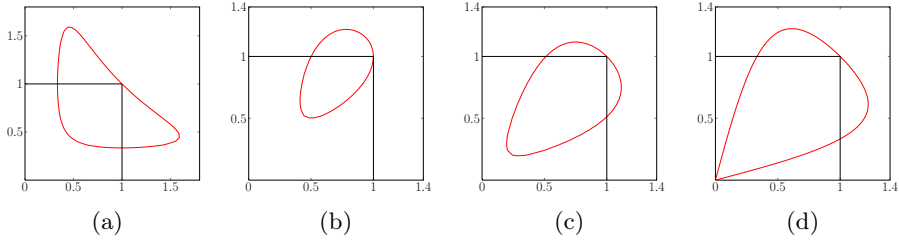}
}
\hfill{}
\caption{Examples of $Q^{+}$. ~\subref{fig:2a} $p_{1,0} = p_{0,1} = \frac{1}{5}$, $p_{-1,-1} = \frac{3}{5}$. ~\subref{fig:2b} $p_{1,0} = \frac{1}{5}$, $p_{0,-1}= p_{-1,1} = \frac{2}{5}$. ~\subref{fig:2c} $p_{1,1} = \frac{1}{62}, p_{-1,1} = p_{1,-1} = \frac{10}{31}, p_{-1,-1} = \frac{21}{62}$. ~\subref{fig:2d} $p_{-1,1} = p_{1,-1} = \frac{1}{4}, p_{-1,-1} = \frac{1}{2}$.}
\label{fig:examplesQ}
\end{figure}

To identify the geometric measures that individually satisfy the balance equations in the interior of the state space, \eqref{eq:interior}, we introduce the polynomial
\begin{equation} \label{eq:Qpoly}
 Q(x,y) = xy\left( \sum_{s = -1}^{1} \sum_{t = -1}^{1} x^{-s} y^{-t} p_{s,t} - 1 \right),
\end{equation}
to capture the notion of balance, \ie $Q(\rho,\sigma)=0$ implies that $m(i,j) = \rho^i\sigma^j, (i,j) \in S$ satisfies~\eqref{eq:interior}. Similarly, we introduce the polynomials $H(x,y)$, $V(x,y)$ to capture the notion of balance for states from the horizontal axis and the vertical axis of the state space, respectively. Let $Q$ denote the set of real $(x,y)$ satisfying $Q(x,y) = 0$, \ie
\begin{equation} \label{eq:Q}
Q = \Big\{ (x,y) \in \mathbb{R}^{2} \mid Q(x,y)=0
\Big\}.
\end{equation}
\noindent Similarly, let $H$ and $V$ denote the set of real $(x,y)$ satisfying $H(x,y) = 0$ and $V(x,y) = 0$ respectively. We are mostly interested in the properties of the algebraic curve $Q$ in $\mathbb{R}_{+}^{2}$, which we denote by $Q^{+}$. Note that $U$ contains those values of $(\rho,\sigma)$ that result in a finite measure. Several examples of $Q^{+}$ are displayed in Figure~\ref{fig:examplesQ}. 


Next, we classify random walks according to~\cite{fayolle1999random}.
\begin{definition}[Singular random walk] \label{def:nonsingular}
Random walk $R$ is called singular if the associated polynomial $Q(x,y)$ is either reducible or of degree $1$ in at least one of the variables.
\end{definition}

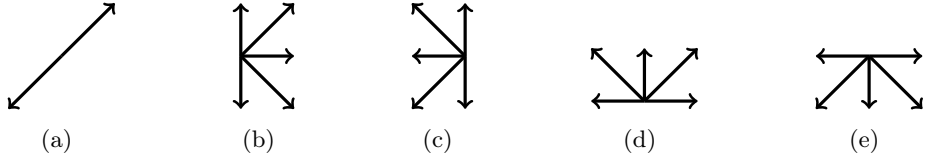
\begin{figure}
\hfill
\subfigure[\label{fig:singular1}]
{
\begin{tikzpicture}[scale=0.7]
\tikzstyle{trans}=[very thick,->]
   \draw[trans] (0,0) to node {$\scriptstyle$} (1,1);
   \draw[trans] (0,0) to node {$\scriptstyle $} (-1,-1);
\end{tikzpicture}
}
\hspace{1cm}
\subfigure[\label{fig:singular2}]
{
\begin{tikzpicture}[scale=0.7]
\tikzstyle{trans}=[very thick,->]
   \draw[trans] (0,0) to node {$\scriptstyle$} (0,-1);
   \draw[trans] (0,0) to node {$\scriptstyle $} (1,-1);
   \draw[trans] (0,0) to node {$\scriptstyle $} (1,0);
   \draw[trans] (0,0) to node {$\scriptstyle $} (1,1);
   \draw[trans] (0,0) to node {$\scriptstyle $} (0,1);
\end{tikzpicture}
}
\hspace{1cm}
\subfigure[\label{fig:singular3}]
{
\begin{tikzpicture}[scale=0.7]
\tikzstyle{trans}=[very thick,->]
   \draw[trans] (0,0) to node {$\scriptstyle$} (0,-1);
   \draw[trans] (0,0) to node {$\scriptstyle $} (-1,-1);
   \draw[trans] (0,0) to node {$\scriptstyle $} (-1,0);
   \draw[trans] (0,0) to node {$\scriptstyle $} (-1,1);
   \draw[trans] (0,0) to node {$\scriptstyle $} (0,1);
\end{tikzpicture}
}
\hspace{1cm}
\subfigure[\label{fig:singular4}]
{
\begin{tikzpicture}[scale=0.7]
\tikzstyle{trans}=[very thick,->]
   \draw[trans] (0,0) to node {$\scriptstyle$} (-1,0);
   \draw[trans] (0,0) to node {$\scriptstyle $} (-1,1);
   \draw[trans] (0,0) to node {$\scriptstyle $} (0,1);
   \draw[trans] (0,0) to node {$\scriptstyle $} (1,1);
   \draw[trans] (0,0) to node {$\scriptstyle $} (1,0);
\end{tikzpicture}
}
\hspace{1cm}
\subfigure[\label{fig:singular5}]
{
\begin{tikzpicture}[scale=0.7]
\tikzstyle{trans}=[very thick,->]
   \draw[trans] (0,0) to node {$\scriptstyle$} (-1,0);
   \draw[trans] (0,0) to node {$\scriptstyle $} (-1,-1);
   \draw[trans] (0,0) to node {$\scriptstyle $} (0,-1);
   \draw[trans] (0,0) to node {$\scriptstyle $} (1,-1);
   \draw[trans] (0,0) to node {$\scriptstyle $} (1,0);
\end{tikzpicture}
}
\hfill{}
\caption{Singular random walks: non-zero transitions in the interior of the state space. \label{fig:singularRW}}
\end{figure}

It was shown in~\cite{fayolle1999random} that a random walk is singular if and only if it has a transition structure that corresponds to one of the cases depicted in Figure~\ref{fig:singularRW}. 

The final piece of notation that we will need deals with the drift of the random walk in the interior of the state space. Let $M_x=\sum_{t = -1}^{1} p_{1,t} - \sum_{t = -1}^{1} p_{-1,t}$ and $M_y=\sum_{s = -1}^{1} p_{s,1} - \sum_{s = -1}^{1} p_{s,-1}$ denote the drift in the horizontal and vertical direction, respectively. In Section~\ref{sec:propertyQ} we will use condition that in an ergodic random walk at least one of the conditions $M_x<0$ or $M_y<0$ holds~\cite{fayolle1999random}. 

Next, we list the assumptions that will be used in the remainder of this paper. A justification for these assumptions will be given below.

\noindent \emph{Assumptions:}
\begin{enumerate}
\item The random walk is non-singular, irreducible, aperiodic and positive recurrent.
\item Measures are positive, finite and absolute convergent in the sense that,
\begin{equation}{\label{eq:absolutelyconvergent}}
\sum_{(\rho, \sigma) \in \Gamma} |\alpha(\rho, \sigma)| \frac{1}{1-\rho} \frac{1}{1-\sigma} < \infty.
\end{equation} 
\item Set $\Gamma$ is a subset of $U = (0,1)^2$ and it satisfies the following technical condition: for any $(\rho, \sigma) \in \Gamma$, there exists a $(w,v) \in \mathbb{N}^2$ such that $\tilde{\rho}^w \tilde{\sigma}^v \neq \rho^w \sigma^v$ for all $(\tilde{\rho}, \tilde{\sigma}) \in \Gamma \backslash (\rho, \sigma)$.
\end{enumerate}

The assumption that random walks are irreducible, aperiodic and positive recurrent ensures that an invariant measure exists. Singular random walks have been analyzed in~\cite{chen2012invariant} for the case that $|\Gamma|<\infty$. It is readily verified that the proofs of~\cite{chen2012invariant} extend to the case that $\Gamma$ has countably infinite cardinality. Therefore, we restrict our attention to non-singular random walks. We are interested in finite measures, in which the sum does not depend on the ordering of the terms. Therefore, we assume absolute convergence~\eqref{eq:absolutelyconvergent}. The assumption that $\Gamma\subset U$ implies that all geometric terms are non-degenerate. The case of degenerate geometric terms was analyzed in~\cite{chen2012invariant} for the case that $|\Gamma|<\infty$. Again, the results from~\cite{chen2012invariant} related to degenerate geometric terms hold for the case that  $\Gamma$ has countably infinite cardinality. Therefore, we assume $\Gamma\subset U$. We will clarify Assumption $3$, \ie the technical condition that we impose on $\Gamma$ in Section~\ref{sec:CandS}.

\section{Algebraic curve $Q$ in $\mathbb{R}^{2}$}{\label{sec:propertyQ}}
In this section we will analyze the algebraic curve $Q$ in $\mathbb{R}^{2}$. Fayolle et al.~\cite{fayolle1999random} have extensively studied the algebraic curve that arises from studying the generating function of the invariant measure and that is defined through $xy(\sum_{s = -1}^{1} \sum_{t = -1}^{1} x^{s} y^{t} \tilde{p}_{s,t} - 1) = 0$. For convenience of notation, we consider the algebraic curve defined by Equation~\eqref{eq:Qpoly} which is the curve considered in~\cite{fayolle1999random} by considering $\tilde p_{i,j} = p_{-i,-j}$. The fact that $\tilde{p}_{i,j}$ might induce a non-ergodic random walk is not a problem, since none of the results from~\cite{fayolle1999random} that will be used in this paper require the random walk to be ergodic.

In this section we will present some of the results from~\cite{fayolle1999random} that will be useful in the sequel as well as a number of new results. The results from~\cite{fayolle1999random} are mostly algebraic of nature. The new results that we present deal with the geometry of $Q$.

The algebraic results that we use from~\cite{fayolle1999random} are expressed in terms of the branch points of the multi-valued algebraic functions $X(y)$ and $Y(x)$ which are defined through
\begin{equation*}
 Q(X(y),y) = Q(x,Y(x)) = 0.
\end{equation*}
These functions are most naturally treated as complex valued functions for complex variables $y$ and $x$, \ie $x, y\in\mathbb{C}$. In particular, the embedding in $\mathbb{C}$ allows us to define the branch points of $X(y)$ and $Y(x)$.

As a first step towards analysis of these branch points, observe that by reordering the terms in $Q(x,y) = 0$ we get
\begin{equation}{\label{eq:Xequation}}
(\sum_{s = -1}^{1} y^{-s + 1} p_{-1,s}) x^2 + (\sum_{s = -1}^{1} y^{-s + 1} p_{0,s} - y) x + (\sum_{s = -1}^{1} y^{-s + 1} p_{1,s}) = 0.
\end{equation}
Therefore, the branch points of $X(y)$ are the roots of $\Delta_x(y) = 0$, where
\begin{equation}{\label{eq:Xdelta}}
\Delta_x(y) = (\sum_{s = -1}^{1} y^{-s + 1} p_{0,s} - y)^2 - 4(\sum_{s = -1}^{1} y^{-s + 1} p_{-1,s})(\sum_{s = -1}^{1} y^{-s + 1} p_{1,s}).
\end{equation}
In similar fashion, by rewriting $Q(x,y) = 0$ into
\begin{equation}{\label{eq:Yequation}}
(\sum_{t = -1}^{1} x^{-t + 1} p_{t,-1}) y^2 + (\sum_{t = -1}^{1} x^{-t + 1} p_{t,0} - x) y + (\sum_{t = -1}^{1} x^{-t + 1} p_{t,1}) = 0,
\end{equation}
it follows that the branch points of $Y(x)$ are the roots of $\Delta_y(x) = 0$,  where
\begin{equation}{\label{eq:Ydelta}}
\Delta_y(x) = (\sum_{t = -1}^{1} x^{-t + 1} p_{t,0} - x)^2 - 4(\sum_{t = -1}^{1} x^{-t + 1} p_{t,-1})(\sum_{t = -1}^{1} x^{-t + 1} p_{t,1}).
\end{equation}

Next, we present two lemmas that fully characterize the location of the branch points of $Y(x)$ and $X(y)$ in terms of the transition probabilities of the random walk. These results provide us with the opportunity to connect the geometry of $Q$ with the interior transition probabilities. The first lemma presented below follows from Lemmas $2.3.8$--$2.3.10$ of~\cite{fayolle1999random}. The result readily follows if one takes into account that in the current paper we consider only ergodic random walks, whereas~\cite{fayolle1999random} also allows for non-ergodic random walks.

\begin{lemma}[Lemmas $2.3.8$--$2.3.10$~\cite{fayolle1999random}] \label{lem:2389}
For all non-singular random walks such that $M_y \neq 0$, $Y(x)$ has four real branch points. Moreover, $Y(x)$ has two branch points $x_1$ and $x_2$ (resp. $x_3$ and $x_4$) inside (resp. outside) the unit circle.

\item For the pair $(x_3, x_4)$, the following classification  holds:
\begin{enumerate}
\item if $p_{-1,0} > 2 \sqrt{p_{-1,-1} p_{-1,1}}$, then $x_3$ and $x_4$ are positive;
\item if $p_{-1,0} = 2 \sqrt{p_{-1,-1} p_{-1,1}}$, then one point is infinite and the other is positive, possibly infinite;
\item if $p_{-1,0} < 2 \sqrt{p_{-1,-1} p_{-1,1}}$, then one point is positive and the other is negative.
\end{enumerate}
\item Similarly, for the pair $(x_1, x_2)$,
\begin{enumerate}
\item if $p_{1,0} > 2 \sqrt{p_{1,-1} p_{1,1}}$, then $x_1$ and $x_2$ are positive;
\item if $p_{1,0} = 2 \sqrt{p_{1,-1} p_{1,1}}$, then one point is $0$ and the second is non-negative;
\item if $p_{1,0} < 2 \sqrt{p_{1,-1} p_{1,1}}$, then one point is positive and the other is negative.
\end{enumerate}

For all non-singular random walks for which $M_y = 0$, one of the branch points of $Y(x)$ is equal to $1$. In addition,
\begin{enumerate}
\item if $M_x < 0$, then two other branch points have a modulus larger than $1$ and the remaining one has a modulus less than $1$;
\item if $M_x > 0$, then two branch points are less than $1$ and the modulus of the remaining one is larger than $1$.
\end{enumerate}
Furthermore, the positivity conditions are the same as the case when $M_y \neq 0$. This lemma holds also for $X(y)$, up to a proper symmetric change of the parameters.
\end{lemma} 
The next lemma deals with multiplicity of the branch points.
\begin{lemma}[Lemma $2.3.10$~\cite{fayolle1999random}]{\label{lem:location}}
The branch points of $X(y)$ and $Y(x)$ with multiplicity $2$ occur only at $0$, $1$ and $\infty$.
\end{lemma}

%

The remainder of this section provides new results on the geometry of $Q$. First, we investigate the possible intersection of $Q$ and $R_{B}$ where $R_B$ is the boundary of the first quadrant.
\begin{lemma}{\label{lem:crossaxis}}
Consider random walk $R$. If $(x,y)\in Q \cap \mathbb{R}_{+}^{2}$ then either $x>0$ and $y>0$ or $x=y=0$, \ie $Q$ cannot cross $R_{B}$ except in the origin.
\end{lemma}
\begin{proof}
If $(x,y)$ is the intersection of $Q$ and $x = 0$, then $y$ must be the root of the following quadratic equation,
\begin{equation}{\label{eq:intersection}}
p_{1,-1} y^2 + p_{1,0} y + p_{1,1} = 0.
\end{equation}
We now show that the roots of~\eqref{eq:intersection} are non-positive by considering all possible choices of $p_{1,-1}, p_{1,0}$ and $p_{1,1}$. If $p_{1,-1} \neq 0$, then~\eqref{eq:intersection} has either no root or two non-positive roots by investigating the relations of the roots using Vieta's formulas. If $p_{1,-1} = 0$ and $p_{1,0} \neq 0$, then~\eqref{eq:intersection} has one non-positive root. If $p_{1,-1} = p_{1,0} = 0$ and $p_{1,1} \neq 0$, then~\eqref{eq:intersection} has no root. The random walk with $p_{1,-1} = p_{1,0} = p_{1,1} = 0$ is excluded because {of the assumption of non-singular random walks}. In similar fashion it follows that $Q$ can only intersect $y = 0$ when $x \leq 0$. Therefore, the only possible intersection of $Q$ and $R_{B}$ is the origin.
\end{proof}

Now we characterize the number of connected components in the first quadrant.
\begin{lemma}{\label{lem:connected}}
Consider random walk $R$. The algebraic curve $Q$ has exactly one closed connected component in $\mathbb{R}_{+}^2$. This component has non-empty intersection with the unit square $U$.
\end{lemma}
\begin{proof}
The algebraic curve $Q$ has exactly one closed connected component in $\mathbb{R}_{+}^2$ follows directly from the fact that $Q(e^x, e^y) < 0$ where $(x,y) \in \mathbb{R}^2$ forms a convex set, as shown in~\cite{latouche2014product,miyazawa2009tail}.

Moreover, the ergodicity conditions for $R$ requires that at least one of the conditions $M_x<0$ or $M_y<0$ holds. Therefore, at least one of the following requirements must be satisfied,
\begin{equation}{\label{eq:intersectionU}}
0 < \frac{\sum_{s = -1}^{1} p_{s,1}}{\sum_{s = -1}^{1} p_{s,-1}} < 1, \quad 0 < \frac{\sum_{t = -1}^{1} p_{1,t}}{\sum_{t = -1}^{1} p_{-1,t}} < 1.
\end{equation}
Since $(1,1) \in Q$, we conclude that this component of $Q$ in the first quadrant has non-empty intersection with $U$.
\end{proof}

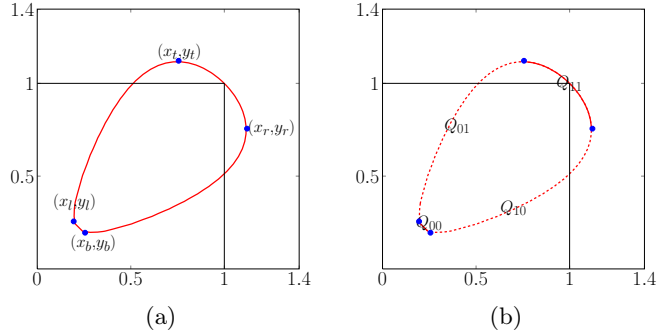
\begin{figure}
\hfill
\subfigure[\label{fig:PPx}]
{
\begin{tikzpicture}[scale = 0.3]
\begin{axis}[%
view={0}{90},
width=4.52083333333333in,
height=4.52083333333333in,
scale only axis,
scale only axis,
xmin=0, xmax=1.4,
ymin=0, ymax=1.4,
xtick = {-0.5,0,0.5,1,1.4},
ytick = {-0.5,0.5,1,1.4}]

\addplot [draw=red, ultra thick] coordinates{ (0.202020202020202,0.353682214701604)(0.198384884176416,0.323232323232323)(0.197077895125914,0.282828282828283)(0.202020202020202,0.246508855441227)(0.203210994526056,0.242424242424242)(0.242424242424242,0.203210994526056)(0.246508855441227,0.202020202020202)(0.282828282828283,0.197077895125914)(0.323232323232323,0.198384884176416)(0.353682214701604,0.202020202020202)(0.363636363636364,0.203895429726532)(0.404040404040404,0.214494599816697)(0.444444444444444,0.225785196573688)(0.484848484848485,0.236997817356168)(0.504404788147808,0.242424242424242)(0.525252525252525,0.250700917059662)(0.565656565656566,0.266180848525009)(0.606060606060606,0.28040464748359)(0.613128657834267,0.282828282828283)(0.646464646464647,0.298405096645817)(0.686868686868687,0.315717377758959)(0.705478529499738,0.323232323232323)(0.727272727272727,0.334879782200734)(0.767676767676768,0.355000115222371)(0.786170830398225,0.363636363636364)(0.808080808080808,0.376913130462884)(0.848484848484849,0.399657695376722)(0.856667825816326,0.404040404040404)(0.888888888888889,0.426152914840109)(0.917591296514448,0.444444444444444)(0.929292929292929,0.453964973403925)(0.969690941847002,0.484848484848485)(0.96969696969697,0.484854384740536)(1.01010101010101,0.522600209342491)(1.01309284479204,0.525252525252525)(1.04839442169988,0.565656565656566)(1.05050505050505,0.568865091941115)(1.07584234204985,0.606060606060606)(1.09090909090909,0.636320777503951)(1.09610259029538,0.646464646464647)(1.10924344886468,0.686868686868687)(1.11585119923513,0.727272727272727)(1.11615563320778,0.767676767676768)(1.11040146355617,0.808080808080808)(1.09884567781367,0.848484848484849)(1.09090909090909,0.867567246365329)(1.08131783270928,0.888888888888889)(1.05778644774189,0.929292929292929)(1.05050505050505,0.939687761385969)(1.02741480688696,0.96969696969697)(1.01010101010101,0.989122644265788)(0.989122644265788,1.01010101010101)(0.96969696969697,1.02741480688696)(0.939687761385969,1.05050505050505)(0.929292929292929,1.05778644774189)(0.888888888888889,1.08131783270928)(0.867567246365329,1.09090909090909)(0.848484848484849,1.09884567781367)(0.808080808080808,1.11040146355617)(0.767676767676768,1.11615563320778)(0.727272727272727,1.11585119923513)(0.686868686868687,1.10924344886468)(0.646464646464647,1.09610259029538)(0.636320777503952,1.09090909090909)(0.606060606060606,1.07584234204985)(0.568865091941115,1.05050505050505)(0.565656565656566,1.04839442169988)(0.525252525252525,1.01309284479204)(0.522600209342491,1.01010101010101)(0.484854384740536,0.96969696969697)(0.484848484848485,0.969690941847002)(0.453964973403925,0.929292929292929)(0.444444444444444,0.917591296514448)(0.426152914840109,0.888888888888889)(0.404040404040404,0.856667825816326)(0.399657695376722,0.848484848484848)(0.376913130462884,0.808080808080808)(0.363636363636364,0.786170830398225)(0.355000115222371,0.767676767676768)(0.334879782200734,0.727272727272727)(0.323232323232323,0.705478529499738)(0.315717377758959,0.686868686868687)(0.298405096645817,0.646464646464647)(0.282828282828283,0.613128657834268)(0.28040464748359,0.606060606060606)(0.266180848525009,0.565656565656566)(0.250700917059662,0.525252525252525)(0.242424242424242,0.504404788147808)(0.236997817356168,0.484848484848485)(0.225785196573688,0.444444444444444)(0.214494599816697,0.404040404040404)(0.203895429726532,0.363636363636364)(0.202020202020202,0.353682214701604)(NaN,NaN)};
\addplot [
color=black,
solid
]
coordinates{
 (0,1)(1,1) 
};

\addplot [
color=black,
solid
]
coordinates{
 (1,0)(1,1) 
};

\node [above] at (axis cs:  0.2, 0.3) {\Huge $\scriptstyle (x_l, y_l)$};
\node [below] at (axis cs:  0.3, 0.2) {\Huge $\scriptstyle (x_b, y_b)$};
\node [above] at (axis cs:  0.76, 1.11) {\Huge $\scriptstyle (x_t, y_t)$};
\node [right] at (axis cs:  1.11, 0.76) {\Huge $\scriptstyle (x_r, y_r)$};
\end{axis}
\filldraw[blue] (1.6,2.1) circle (3pt) (2.1,1.6) circle (3pt) (9.2,6.2) circle (3pt) (6.2,9.2) circle (3pt);
\end{tikzpicture}
}
\subfigure[\label{fig:PPc}]
{
\begin{tikzpicture}[scale = 0.3]
\begin{axis}[%
view={0}{90},
width=4.52083333333333in,
height=4.52083333333333in,
scale only axis,
scale only axis,
xmin=0, xmax=1.4,
ymin=0, ymax=1.4,
xtick = {-0.5,0,0.5,1,1.4},
ytick = {-0.5,0.5,1,1.4}]

\addplot [draw=red, ultra thick, dashed] coordinates{ (0.202020202020202,0.353682214701604)(0.198384884176416,0.323232323232323)(0.197077895125914,0.282828282828283)(0.202020202020202,0.246508855441227)(0.203210994526056,0.242424242424242)(0.242424242424242,0.203210994526056)(0.246508855441227,0.202020202020202)(0.282828282828283,0.197077895125914)(0.323232323232323,0.198384884176416)(0.353682214701604,0.202020202020202)(0.363636363636364,0.203895429726532)(0.404040404040404,0.214494599816697)(0.444444444444444,0.225785196573688)(0.484848484848485,0.236997817356168)(0.504404788147808,0.242424242424242)(0.525252525252525,0.250700917059662)(0.565656565656566,0.266180848525009)(0.606060606060606,0.28040464748359)(0.613128657834267,0.282828282828283)(0.646464646464647,0.298405096645817)(0.686868686868687,0.315717377758959)(0.705478529499738,0.323232323232323)(0.727272727272727,0.334879782200734)(0.767676767676768,0.355000115222371)(0.786170830398225,0.363636363636364)(0.808080808080808,0.376913130462884)(0.848484848484849,0.399657695376722)(0.856667825816326,0.404040404040404)(0.888888888888889,0.426152914840109)(0.917591296514448,0.444444444444444)(0.929292929292929,0.453964973403925)(0.969690941847002,0.484848484848485)(0.96969696969697,0.484854384740536)(1.01010101010101,0.522600209342491)(1.01309284479204,0.525252525252525)(1.04839442169988,0.565656565656566)(1.05050505050505,0.568865091941115)(1.07584234204985,0.606060606060606)(1.09090909090909,0.636320777503951)(1.09610259029538,0.646464646464647)(1.10924344886468,0.686868686868687)(1.11585119923513,0.727272727272727)(1.11615563320778,0.767676767676768)(1.11040146355617,0.808080808080808)(1.09884567781367,0.848484848484849)(1.09090909090909,0.867567246365329)(1.08131783270928,0.888888888888889)(1.05778644774189,0.929292929292929)(1.05050505050505,0.939687761385969)(1.02741480688696,0.96969696969697)(1.01010101010101,0.989122644265788)(0.989122644265788,1.01010101010101)(0.96969696969697,1.02741480688696)(0.939687761385969,1.05050505050505)(0.929292929292929,1.05778644774189)(0.888888888888889,1.08131783270928)(0.867567246365329,1.09090909090909)(0.848484848484849,1.09884567781367)(0.808080808080808,1.11040146355617)(0.767676767676768,1.11615563320778)(0.727272727272727,1.11585119923513)(0.686868686868687,1.10924344886468)(0.646464646464647,1.09610259029538)(0.636320777503952,1.09090909090909)(0.606060606060606,1.07584234204985)(0.568865091941115,1.05050505050505)(0.565656565656566,1.04839442169988)(0.525252525252525,1.01309284479204)(0.522600209342491,1.01010101010101)(0.484854384740536,0.96969696969697)(0.484848484848485,0.969690941847002)(0.453964973403925,0.929292929292929)(0.444444444444444,0.917591296514448)(0.426152914840109,0.888888888888889)(0.404040404040404,0.856667825816326)(0.399657695376722,0.848484848484848)(0.376913130462884,0.808080808080808)(0.363636363636364,0.786170830398225)(0.355000115222371,0.767676767676768)(0.334879782200734,0.727272727272727)(0.323232323232323,0.705478529499738)(0.315717377758959,0.686868686868687)(0.298405096645817,0.646464646464647)(0.282828282828283,0.613128657834268)(0.28040464748359,0.606060606060606)(0.266180848525009,0.565656565656566)(0.250700917059662,0.525252525252525)(0.242424242424242,0.504404788147808)(0.236997817356168,0.484848484848485)(0.225785196573688,0.444444444444444)(0.214494599816697,0.404040404040404)(0.203895429726532,0.363636363636364)(0.202020202020202,0.353682214701604)(NaN,NaN)};

\addplot [draw=red, ultra thick] coordinates{ (0.203210994526056,0.242424242424242)(0.242424242424242,0.203210994526056)};

\addplot [draw=red, ultra thick] coordinates{ (1.11615563320778,0.767676767676768)(1.11040146355617,0.808080808080808)(1.09884567781367,0.848484848484849)(1.09090909090909,0.867567246365329)(1.08131783270928,0.888888888888889)(1.05778644774189,0.929292929292929)(1.05050505050505,0.939687761385969)(1.02741480688696,0.96969696969697)(1.01010101010101,0.989122644265788)(0.989122644265788,1.01010101010101)(0.96969696969697,1.02741480688696)(0.939687761385969,1.05050505050505)(0.929292929292929,1.05778644774189)(0.888888888888889,1.08131783270928)(0.867567246365329,1.09090909090909)(0.848484848484849,1.09884567781367)(0.808080808080808,1.11040146355617)(0.767676767676768,1.11615563320778)};

\addplot [
color=black,
solid
]
coordinates{
 (0,1)(1,1) 
};

\addplot [
color=black,
solid
]
coordinates{
 (1,0)(1,1) 
};
\node at (axis cs:  0.25,  0.25) {\Huge $\scriptstyle Q_{00}$};
\node at (axis cs:  0.7,  0.33) {\Huge $\scriptstyle Q_{10}$};
\node at (axis cs:  1,  1) {\Huge $\scriptstyle Q_{11}$};
\node at (axis cs:  0.4,  0.77) {\Huge $\scriptstyle Q_{01}$};
\end{axis}
\filldraw[blue] (1.6,2.1) circle (3pt) (2.1,1.6) circle (3pt) (9.2,6.2) circle (3pt) (6.2,9.2) circle (3pt);
\end{tikzpicture}
} 
\hfill{}
\caption{$Q^{+}$ for the random walk from Figure~\ref{fig:2c}.~\subref{fig:PPx} Branch points of $Q^{+}$.~\subref{fig:PPc} Partition of $Q^{+}$.
\label{fig:PPQ}}

\end{figure}

Now we know that $Q^{+}$, the intersection of $Q$ and $\mathbb{R}_{+}^{2}$, is a connected component. Denote the branch points of $Y(x)$ and $X(y)$ on $Q^{+}$ by $x_l, x_r$ with $x_l < x_r$ and $y_b, y_t$ with $y_b < y_t$ respectively. Let $y_l, y_r, x_b, x_t$ satisfy $(x_l, y_l), (x_r, y_r), (x_t, y_t), (x_b, y_b) \in Q$, see Figure~\ref{fig:PPx}. We will refer to $(x_l, y_l), (x_r, y_r), (x_t, y_t), (x_b, y_b)$ as branch points of $Q^+$. From Lemma~\ref{lem:2389}, we know that $0 \leq x_l \leq 1  \leq x_r$, $0 \leq y_b \leq 1 \leq y_t$. Since we are only interested in finite measures, we only consider $Q^{+}$ in $\bar{U} = [0,1)^2$. Lemma~\ref{lem:connected} states that $Q^{+}_U = Q^{+} \cap U$ is a non-empty set for an ergodic random walk with non-zero drift. We first start the analysis of $Q^{+}$.
\begin{definition}[Partition of $Q^{+}$]
The partition $\{Q_{00}, Q_{01}, Q_{10}, Q_{11}\}$ of $Q^{+}$ is defined as follows: $Q_{00}$ is the part of $Q$ connecting $(x_l, y_l)$ and $(x_b, y_b)$; $Q_{10}$ is the part of $Q$ connecting $(x_b, y_b)$ and $(x_r, y_r)$; $Q_{01}$ is the part of $Q$ connecting $(x_l, y_l)$ and $(x_t, y_t)$; $Q_{11}$ is the part of $Q$ connecting $(x_r, y_r)$ and $(x_t, y_t)$.
\end{definition}
An example of the partition of $Q^{+}$ is given in Figure~\ref{fig:PPc}.
The monotonicity of $X(y)$ and $Y(x)$ will play a crucial role in analyzing the structure of $\Gamma$, we state it here without proof because it immediately follows from the fact that  $\{(x,y) \in \mathbb{R}^2 : Q(e^x, e^y) < 0\}$ is convex. 
\begin{lemma}{\label{lem:monotonicity}}
Consider $(x,y) \in Q_{i,1-i}$ and $(\tilde{x}, \tilde{y}) \in Q_{i,1-i}$ where $i = 0,1$, if $\tilde{x} > x$, then $\tilde{y} > y$. Consider $(x,y) \in Q_{i,i}$ and $(\tilde{x}, \tilde{y}) \in Q_{i,i}$ where $i = 0,1$, if $\tilde{x} > x$, then $\tilde{y} < y$.
\end{lemma}

Next, we turn our attention to the singularities of $Q$. 
\begin{definition}[Singularity of $Q$]\label{def:singularity}
Point $(x, y)\in Q$ is a singularity of multiplicity $m$, $m>1$, iff at $(x, y)$ all partial derivatives of $Q(x,y)$ of order less than $m$ vanish and at least one partial derivative of order $m$ is non-zero.
\end{definition}

We will see below that if one set from the partition of $Q^{+}$ is empty, then the curve $Q^{+}$ will have a singularity. The singularity plays an important role in the analysis later.
\begin{figure}
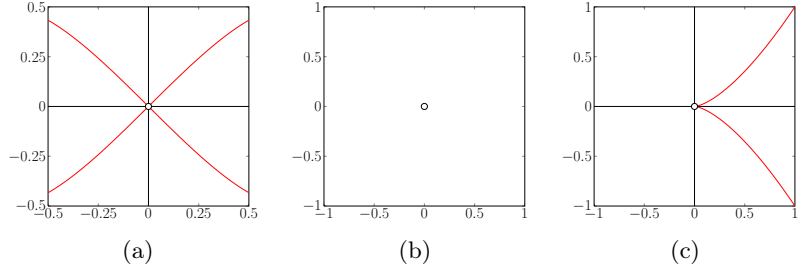

\hfill
\subfigure[\label{fig:DPa}]
{
\input{crunode.tikz}
} 
\subfigure[\label{fig:DPb}]
{
%
%
%
%
\begin{tikzpicture}[scale = 0.23]

\begin{axis}[%
view={0}{90},
width=4.52083333333333in,
height=4.52083333333333in,
scale only axis,
xmin=-1, xmax=1,
ymin=-1, ymax=1,
xtick = {-1,-0.5,0,0.5,1},
ytick = {-1,-0.5,0,0.5,1}]
\addplot [
color=blue,
mark size=5.0pt,
only marks,
mark=*,
mark options={solid,fill=white,draw=black}
]
coordinates{
 (0,0) 
};

\end{axis}
\end{tikzpicture}
}
\subfigure[\label{fig:DPc}] 
{
\input{cusp.tikz}
}
\hfill{}
\caption{Types of real double point. ~\subref{fig:DPa} crunode. ~\subref{fig:DPb} acnode. ~\subref{fig:DPc} ordinary cusp. \label{fig:doublepoints}}
\end{figure} 
\begin{lemma}{\label{lem:BPandS}}
For all random walks with non-zero drift, $(x,y)$ is a singularity of $Q^{+}$ iff it is a crunode of order 2 and $x$ and $y$ are branch points of multiplicity $2$ of $Y(x)$ and $X(y)$ respectively.
\end{lemma}
\begin{proof}
We prove by contradiction that it is not possible to have a singularity of order larger than $2$. Suppose that $(\tilde x,\tilde y)\in Q^+$ is a singularity of order larger than $2$. From Lemma~\ref{lem:crossaxis} it follows that we need to consider the cases i) $\tilde x>0$ and $\tilde y>0$, and ii) $(\tilde x,\tilde y)=(0,0)$.  If $\tilde x>0$ and $\tilde y>0$ it follows from 
\begin{equation*}
\frac{\partial^2 Q(x,y)}{\partial^2 x} = \sum_{t = -1}^{1} p_{-1,t} y^{-t + 1} = 0, \quad\text{and}\quad
\frac{\partial^2 Q(x,y)}{\partial^2 y} = \sum_{s = -1}^{1} p_{s,-1} x^{-s + 1} = 0,
\end{equation*}  
that $p_{-1,1} = p_{-1,0} = p_{-1,-1} = p_{0,-1} = p_{1,-1} = 0$, which leads to a non-ergodic random walk. For $(\tilde x,\tilde y)=(0,0)$ it follows from
\begin{equation*}
\frac{\partial^2 Q(x,y)}{\partial x\partial y} = 4xy p_{-1,-1} + 2x p_{-1,0} + 2y p_{0,-1} + p_{0,0} - 1 = 0,
\end{equation*} 
that $p_{00} = 1$, which leads to a random walk that is not irreducible. This concludes the proof that a singularity has at most order $2$. 

Next, we demonstrate that if $(x, y)$ is a singularity, then $x$ and $y$ are branch points of $Y(x)$ and $X(y)$ respectively. By combining $Q(x,y)=0$ with
\begin{equation*}
\frac{\partial Q(x,y)}{\partial x} = 2x (\sum_{t = -1}^{1} p_{-1,t} y^{-t + 1}) + (\sum_{t = -1}^{1} p_{0,t} y^{-t + 1} - y) = 0\label{eq:firstx}
\end{equation*}
we obtain
\begin{equation*}\label{eq:singularbranch}
\sum_{t = -1}^{1} p_{-1,t} y^{-t + 1} x^2 = \sum_{t = -1}^{1} p_{1,t} y^{-t+1},
\end{equation*}
which means $x$ is the root of $\Delta_y(x) = 0$, defined in~\eqref{eq:Ydelta} and therefore a branch point of $Y(x)$. Similarly, it follows from $Q(x,y)=0$ and $\partial Q(x,y)/\partial y=0$ that $y$ is a branch point of $X(y)$.

Now, we are ready to prove that a singularity $(x,y)$ is a crunode. For more information on the classification of singularities of algebraic curves, see \eg~\cite{gibson1998elementary}. An illustration of all possible singularities of order $2$ is given in Figure~\ref{fig:doublepoints}. Note, that the figure does not include a ramphoid cusp, since it has order larger than $2$. A singularity cannot be an ordinary cusp, because $x$ and $y$ are branch points of $Y(x)$ and $X(y)$ respectively.  Moreover, $(x,y)$ is not an acnode because $Q^{+}_U$ is non-empty due to Lemma~\ref{lem:connected}. Therefore, a singularity is a crunode.

The final result in this lemma follows from the observation that if $x$ and $y$ are branch points of $Y(x)$ and $X(y)$ respectively and $(x,y)$ is a crunode then $x$ and $y$ must have multiplicity two.
\end{proof}

\begin{theorem}{\label{thm:singularity8}}
The algebraic curve $Q$ has a singularity in $\bar{U}$ if and only if $p_{0,1} = p_{1,1} = p_{1,0} = 0$, in which case this singularity is located in the origin.
\end{theorem}
\begin{proof}
Lemma~\ref{lem:BPandS} states that $(x,y)$ is a singularity of $Q^{+}$ if and only if it is a crunode of order 2 and $x$ and $y$ are branch points of multiplicity $2$ of $Y(x)$ and $X(y)$ respectively. Therefore, we only need to consider $(x,y)$ where $x$ and $y$ are the multiple roots of $\Delta_y(x) = 0$ and $\Delta_x(y) = 0$ respectively. A multiple root of $\Delta_y(x) = 0$ and $\Delta_x(y) = 0$ can only occur at $x = 0, 1$ or $\infty$ and $y = 0, 1$ or $\infty$, respectively, due to Lemma~\ref{lem:location}. Therefore, $x = 0$ and $y = 0$ must be multiple roots of $\Delta_y(x) = 0$ and $\Delta_x(y) = 0$, respectively, if there is a singularity in $\bar{U}$. We know from~\cite[Lemma $2.3.10$]{fayolle1999random} that $\Delta_y(x) = 0$ has a multiple root at $0$ if and only if one of the following holds:
\begin{equation}{\label{eq:conditionx1}}
p_{-1,0} = p_{-1,1} = p_{0,1} = 0,
\end{equation}
\begin{equation}{\label{eq:conditionx2}}
p_{1,0} = p_{1,1} = p_{0,1} = 0,
\end{equation}
\begin{equation}{\label{eq:conditionx3}}
p_{-1,-1} = p_{0,-1} = p_{1,-1} = 0,
\end{equation}
and $\Delta_x(y) = 0$ has a multiple root at $0$ if and only if one of the following holds:
\begin{equation}{\label{eq:conditiony1}}
p_{0,-1} = p_{1,-1} = p_{1,0} = 0,
\end{equation}
\begin{equation}{\label{eq:conditiony2}}
p_{0,1} = p_{1,1} = p_{1,0} = 0,
\end{equation}
\begin{equation}{\label{eq:conditiony3}}
p_{-1,-1} = p_{-1,0} = p_{-1,1} = 0.
\end{equation}
Conditions~\eqref{eq:conditionx3} and~\eqref{eq:conditiony3} lead to a singular random walk. The combinations of conditions~\eqref{eq:conditionx1} and~\eqref{eq:conditiony1},~\eqref{eq:conditionx1} and~\eqref{eq:conditiony2},~\eqref{eq:conditionx2} and~\eqref{eq:conditiony1} will lead to singular random walks as well. Since by Assumption 1 the random walk is non-singular, the algebraic curve $Q$ has a singularity in $\bar{U}$ if and only if $p_{0,1} = p_{1,1} = p_{1,0} = 0$, in which case it is located in the origin.
\end{proof}

\section{Constraints on invariant measures and random walks}{\label{sec:CandS}}
In this section, we will first demonstrate that $\Gamma\subset Q^+_U$, \ie candidate geometric measures individually satisfy the balance equations in the interior of the state space. Next we characterize the structure of $\Gamma$ that may lead to an invariant measure. Then we will provide necessary conditions on the transition probabilities of a random walk to allow for an infinite sum of  geometric terms to constitute the invariant measure. Finally, we demonstrate that it is necessary to have at least one negative coefficient in an invariant measure that is an infinite sum of geometric terms.

The results in Subsections~\ref{ssec:gammainQ} and~\ref{ssec:gammapairwise} make use of the following result that is a special case of a result from~\cite{brown1960absolutely}.
\begin{theorem}[Theorem $1$, Lemma $1$~\cite{brown1960absolutely}]{\label{thm:brown}}
Let $d$ be a positive integer. Consider a real measure $\mu:\mathbb{R}^d\to\mathbb{R}$ with bounded compact support $K$. If 
\begin{equation}{\label{eq:intpzero}}
\int_{\mathbb{R}^d} P(x) \,d \mu(x) = 0
\end{equation}
for all polynomials $P$, then $\mu = 0$.
\end{theorem}
\subsection{Constraints on set $\Gamma$} \label{ssec:gammainQ}
We demonstrate in this subsection that only $\Gamma\subset Q^{+}_U$ may induce an invariant measure for a random walk. In particular we prove that this is true under the technical condition that for any $(\rho, \sigma) \in \Gamma$, there exists a $(w,v) \in \mathbb{N}^2$ such that $\tilde{\rho}^w \tilde{\sigma}^v \neq \rho^w \sigma^v$ for all $(\tilde{\rho}, \tilde{\sigma}) \in \Gamma \backslash (\rho, \sigma)$, which was already introduced as Assumption 3 in Section~\ref{sec:model}. The interpretation of this condition is that each of the elements from $\Gamma$ can in a sense by separated from all other elements in $\Gamma$ by selecting the proper $(v,w)\in\mathbb{N}^2$. This separating principle will be used in our proofs to demonstrate that each term $(\rho,\sigma)$ in $Q$ must individually satisfy the balance equations in the interior of the state space, hence $\Gamma \subset Q^{+}_U$.     

Next, we present the main result of this subsection.
 
\begin{theorem}{\label{thm:productterm}}
If the invariant measure for a random walk in the quarter-plane is induced by $\Gamma$, where for any $(\rho, \sigma) \in \Gamma$ there exists a $(w,v) \in \mathbb{N}^2$ such that $\tilde{\rho}^w \tilde{\sigma}^v \neq \rho^w \sigma^v$ for all $(\tilde{\rho}, \tilde{\sigma}) \in \Gamma \backslash (\rho, \sigma)$, then $\Gamma \subset Q_U^{+}$.
\end{theorem}

\begin{proof}
Consider $(\rho, \sigma) \in \Gamma$ such that there exists a $(w,v) \in \mathbb{N}_{+}^2$ such that $\tilde{\rho}^w \tilde{\sigma}^v \neq \rho^w \sigma^v$ for any $(\tilde{\rho}, \tilde{\sigma})\in \Gamma \backslash (\rho, \sigma)$. We now partition $\Gamma \backslash \{(\rho, \sigma)\}$ into elements $\Gamma_1$, $\Gamma_2, \cdots, \Gamma_k, \cdots$ as follows. If $\rho_m^w \sigma_m^v = \rho_n^w \sigma_n^v$, then $(\rho_m, \sigma_m)$ and $(\rho_n, \sigma_n)$ will be put into the same element in the partition. Moreover, we arbitrarily choose one geometric term from each element as the representative, which is denoted by $(\rho(\Gamma_k), \sigma(\Gamma_k))$ where $k = 1,2,3, \cdots$. For convenience of notation, we denote $\Gamma_0 = (\rho, \sigma)$.
 
Since $m$ satisfies balance equation~\eqref{eq:interior} in the interior of the state space, \ie $(i,j)$ with $i \geq 0$ and $j \geq 0$,
\begin{equation*}{\label{eq:generalbalance}}
\sum_{(\rho, \sigma) \in \Gamma} \rho^i \sigma^j [\alpha(\rho, \sigma) \rho \sigma (1 - \sum_{s = -1}^{1} \sum_{t = -1}^1 \rho^{-s} \sigma^{-t} p_{s,t})] = 0.
\end{equation*}

We now consider the balance equation for states $(dw + 1, dv + 1)$ where $d = 1,2,3, \cdots$,

\begin{equation}{\label{eq:Ginfinite}}
\sum_{k = 0}^{\infty} [\rho(\Gamma_k)^w \sigma(\Gamma_k)^v]^d [\sum_{(\rho, \sigma) \in \Gamma_k}\alpha(\rho, \sigma) (\rho \sigma - \sum_{s = -1}^{1} \sum_{t = -1}^1 \rho^{1-s} \sigma^{1-t} p_{s,t})]= 0.
\end{equation}
Below we will apply Theorem~\ref{thm:brown}. For this we require absolutely convergence of the series
\begin{equation*}
\sum_{(\rho, \sigma) \in \Gamma_k} \alpha(\rho, \sigma) (\rho \sigma - \sum_{s = -1}^{1} \sum_{t = -1}^1 \rho^{1-s} \sigma^{1-t} p_{s,t}),
\end{equation*}
where $k = 0,1,2,\cdots$, which we show first. Because of assumption~\eqref{eq:absolutelyconvergent}, \ie the absolute convergence of the terms of which the sum is $m$, we have
\begin{align*}
&\sum_{(\rho, \sigma) \in \Gamma_k} \alpha(\rho, \sigma) (\rho \sigma - \sum_{s = -1}^{1} \sum_{t = -1}^1 \rho^{1-s} \sigma^{1-t} p_{s,t})  \\
< & B \sum_{(\rho, \sigma) \in \Gamma}|\alpha(\rho, \sigma)| \frac{1}{1 - \rho} \frac{1}{1 - \sigma}\\
< & \infty,
\end{align*}
where $B$ is a finite positive constant. We define a real measure $\mu$ on $\mathbb{R}^2$ as
\begin{equation*}
\mu(\rho, \sigma) =
\begin{cases}
  \sum_{(\rho, \sigma) \in \Gamma_k}\alpha(\rho, \sigma)(\rho \sigma - \sum_{s = -1}^{1} \sum_{t = -1}^1 \rho^{1-s} \sigma^{1-t} p_{s,t})\\
  \quad \text{if} \quad (\rho,\sigma) = (\rho(\Gamma_k), \sigma(\Gamma_k)) \quad \text{for} \quad k = 0,1,2, \cdots \\
  0 \quad \text{otherwise}.
\end{cases}
\end{equation*}
We can now write~\eqref{eq:Ginfinite} as
\begin{align*}
&\sum_{k = 0}^{\infty} [\rho(\Gamma^k)^{w} \sigma(\Gamma^k)^{v}]^d \mu(\rho, \sigma) = \int [\rho(\Gamma^k)^{w} \sigma(\Gamma^k)^{v}]^d \,d \mu(\rho, \sigma) = 0,
\end{align*}
for $d = 0,1,2 \cdots$. This indicates that
\begin{equation*}
\int P(\rho, \sigma) \,d \mu(\rho, \sigma) = 0
\end{equation*}
for all $P(\rho, \sigma) = (\rho^w \sigma^v)^d$ where $d =0, 1, 2 \cdots$. Hence $\int P(\rho, \sigma) \,d \mu(\rho, \sigma) = 0$ for all polynomials. Moreover, for $k =0, 1, 2, \cdots$, the following series is absolutely convergent, $\sum_{(\rho, \sigma) \in \Gamma^k}\alpha(\rho, \sigma)(\rho \sigma - \sum_{s = -1}^{1} \sum_{t = -1}^1 \rho^{1-s} \sigma^{1-t} p_{s,t})$. Hence, the compact support of this sequence is a bounded interval. Therefore, by using Theorem~\ref{thm:brown}, $\mu = 0$, thus $1 - \sum_{s = -1}^{1} \sum_{t = -1}^1 \rho^{-s} \sigma^{-t} p_{s,t} = 0$ because of the assumption of non-degenerate geometric terms, \ie $\rho \sigma \neq 0$, which completes the proof.
\end{proof}

\subsection{Structure of $\Gamma$} \label{ssec:gammapairwise}
In this section, we consider the structure of $\Gamma$. The proofs in this and subsequent sections are based on the notion of uncoupled partitions, which is introduced first.
\begin{definition}[Uncoupled partition]
A partition $\{\Gamma_1, \Gamma_2, \cdots\}$ of $\Gamma$ is \emph{horizontally uncoupled} if $(\rho, \sigma) \in \Gamma_p$ and $(\tilde{\rho}, \tilde{\sigma}) \in \Gamma_q$ for $p \neq q$, implies that $\tilde{\rho} \neq \rho$, \emph{vertically uncoupled} if $(\rho, \sigma) \in \Gamma_p$ and $(\tilde{\rho}, \tilde{\sigma}) \in \Gamma_q$ for $p \neq q$, implies $\tilde{\sigma} \neq \sigma$, and \emph{uncoupled} if it is both horizontally and vertically uncoupled.
\end{definition}

\begin{figure}
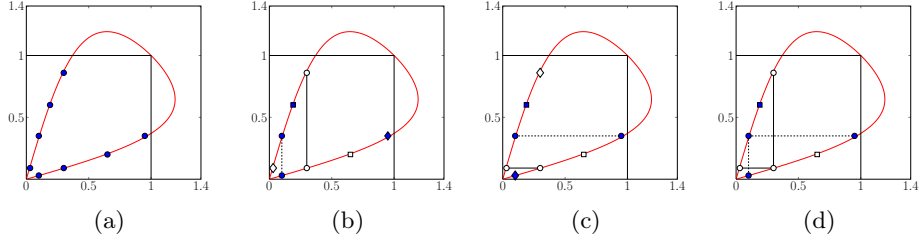

\hfill
\subfigure[
{\label{fig:CandGamma}}
]
{
\input{partition1.tikz}
} 
\subfigure[
 {\label{fig:HUP}}
 ] 
 {
 \input{partition3.tikz}
 }
\subfigure[
{\label{fig:VUP}}
]
  {
  \input{partition4.tikz}
  }
\subfigure[
{\label{fig:MUP}}
]
{
\input{partition2.tikz}
}
\hfill{}
\caption{Partitions of set $\Gamma$.~\subref{fig:CandGamma} curve $Q^{+}$ of Figure~\ref{fig:2d} and $\Gamma \subset Q^{+}$ as dots.~\subref{fig:HUP} horizontally uncoupled partition with $6$ sets.~\subref{fig:VUP} vertically uncoupled partition with $6$ sets.~\subref{fig:MUP} uncoupled partition with $4$ sets. Different sets are marked by different symbols.}
\label{fig:allfigures}
\end{figure} 
We call a partition with the largest number of sets a maximal partition. 
\begin{lemma}[Lemma $1$~\cite{chen2012invariant}]{\label{lem:unique}}
The maximal vertically uncoupled partition, the maximal horizontally uncoupled partition and the maximal uncoupled partition are unique.
\end{lemma}
Examples of a maximal horizontally uncoupled partition, of a maximal vertically uncoupled partition and of a maximal uncoupled partition can be found in Figure~\ref{fig:allfigures}. Let $H[\Gamma]$ denote the number of elements in the maximal horizontally uncoupled partition for set $\Gamma$, and let $\Gamma_p^h$, $p=1,\dots,H[\Gamma]$ denote the sets themselves, where elements of $\Gamma_p^h$ have common horizontal coordinate $\rho(\Gamma_p^h)$. The maximal vertically uncoupled partition of set $\Gamma$ has $V[\Gamma]$ sets, $\Gamma_q^v$, $q=1,\cdots,V[\Gamma]$, where elements of $\Gamma_q^v$ have common vertical coordinate $\sigma(\Gamma_q^v)$. The maximal uncoupled partition of set $\Gamma$ is denoted by $\{\Gamma_k^g\}_{k=1}^{G[\Gamma]}$. The $H[\Gamma], V[\Gamma]$ and $G[\Gamma]$ are allowed to be infinite.

We now make two observations on the structure of an element $\Gamma_k^g$ from the maximal uncoupled partition. Firstly, for any $(\rho,\sigma)\in\Gamma_k^g$ there always exist either $(\rho,\tilde\sigma)\in\Gamma_k^g$, with $\tilde{\sigma} \neq \sigma$ or $(\tilde\rho,\sigma)\in\Gamma_k^g$ with $\tilde{\rho} \neq \rho$. Secondly, the degree of $Q(\rho,\sigma)$ is at most two in each variable. This means, for instance, that if $(\rho,\sigma)\in\Gamma_k^g$ and $(\rho,\tilde\sigma)\in\Gamma_k^g$, $\tilde\sigma\neq\sigma$, then there does not exist $(\rho,\hat\sigma)\in\Gamma_k^g$, where $\hat{\sigma} \neq \sigma$ and $\hat{\sigma} \neq \tilde{\sigma}$. By repeating the above two arguments for other elements in $\Gamma_k^g$ it follows that $\Gamma_k^g$ must have a pairwise-coupled structure. An example of such a set is $\Gamma_k^g = \{(\rho_k, \sigma_k), k = 1,2,3 \cdots \}$, where
\begin{equation} \label{eq:pairwisexample}
 \rho_1 = \rho_2, \sigma_1 > \sigma_2, \rho_2 > \rho_3, \sigma_2 = \sigma_3, \rho_3 = \rho_4, \sigma_3 > \sigma_4, \cdots.
\end{equation}
The above discussion leads to the definition of a pairwise-coupled set in terms of the number of sets in a maximal uncoupled partition.
\begin{definition}[Pairwise-coupled set]
A set $\Gamma \subset Q^{+}_{U}$ is pairwise-coupled if and only if the maximal uncoupled partition of $\Gamma$ contains only one set.
\end{definition}
Note that the above definition implies that each of the sets in $\{\Gamma_k^g\}_{k=1}^{G[\Gamma]}$ is pairwise coupled.

We are now ready to show that if the measure induced by $\Gamma$ is the invariant measure, then $\Gamma$ must be the union of finitely many pairwise-coupled sets each with countably infinite cardinality. We first introduce some additional notation. For any set $\Gamma_p^h$ from the maximal horizontally uncoupled partition of $\Gamma$, let
\begin{align*}
B^h(\Gamma_p^h) &= \sum_{(\rho, \sigma)\in\Gamma_p^h}\alpha(\rho, \sigma)[\sum_{s=-1}^1 \big(\rho^{1-s} h_s+\rho^{1-s}\sigma p_{s,-1}\big) - \rho].  
\end{align*}
For any set $\Gamma_q^v$ from the maximal vertically uncoupled partition of $\Gamma$, let
\begin{align*}
B^v(\Gamma_q^v) &= \sum_{(\rho, \sigma)\in\Gamma_q^v}\alpha(\rho, \sigma)[\sum_{t=-1}^1 \big(\sigma^{1-t} v_t+\rho \sigma^{1-t} p_{-1,t}\big) - \sigma]. 
\end{align*}

\begin{lemma} \label{lem:MUPIM}
Consider the random walk $R$, if the invariant measure $m$ of $R$ is $m(i,j) = \sum_{(\rho, \sigma) \in \Gamma} \alpha(\rho, \sigma)\rho^i \sigma^j$, then we have  $B^h(\Gamma_p^h) = 0$ and $B^v(\Gamma_q^v) = 0$ for all $p=1,\dots,H[\Gamma]$ and $q=1,\dots,V[\Gamma]$, respectively.
\end{lemma}
\begin{proof}
Since $m$ is the invariant measure of $R$, $m$ satisfies the balance equations at state $(i,0)$ for $i = 1,2,3 \cdots$. Therefore,
\begin{align} 
0
&= \sum_{s=-1}^1 \big[m(i-s,0)h_s + m(i-s,1)p_{s,-1}\big] - m(i,0) \notag \\
&= \sum_{(\rho, \sigma)\in\Gamma} \alpha(\rho,\sigma) [\sum_{s=-1}^1 \big(\rho^{i-s} h_s+\rho^{i-s}\sigma p_{s,-1}\big) - \rho^i ] \notag  \\
&= \sum_{p=1}^{\infty} \rho(\Gamma_p^h)^{i}\sum_{(\rho,\sigma) \in\Gamma_p^h}\alpha(\rho,\sigma)[\sum_{s=-1}^1 \big(\rho^{-s} h_s+\rho^{-s}\sigma p_{s,-1}\big) - 1] \notag \\
&= \sum_{p=1}^{\infty} \rho(\Gamma_p^h)^{i-1} B^h(\Gamma_p^h). \label{eq:A}
\end{align}
We now show the absolute convergence of sequence $\{B^h (\Gamma_p^h)\}$, $p = 1, \cdots, H[\Gamma]$. Because of the assumption that $m(i,j)$ is a finite measure, we have
\begin{align}
&\sum_{p = 1}^{\infty} |B^h (\Gamma_p^h)| \notag \\
\leq & \sum_{k = 1}^{\infty} \alpha(\rho_k, \sigma_k)[\sum_{s=-1}^1 \big(\rho_k^{1-s} h_s+\rho_k^{1-s}\sigma p_{s,-1}\big) - \rho_k] \notag \\
< & \infty, \notag
\end{align}
the last inequality holds due to the same reasoning that is in the proof of Theorem~\ref{thm:productterm}. If we define a real measure $\mu$ as
\begin{equation*}
\mu(\rho) =
\begin{cases}
 B^h(\Gamma_p^h) \quad \text{if} \quad \rho = \rho(\Gamma_p^h), \\
  0 \quad \text{otherwise}.
\end{cases}
\end{equation*}
We can now write~\eqref{eq:A} as
\begin{equation*}
\sum_{p = 1}^{\infty} \rho(\Gamma_p^h)^{i - 1} \mu(\rho(\Gamma_p^h)) = \int \rho(\Gamma_p^h)^{i - 1} \,d \mu(\rho(\Gamma_p^h)) = 0,
\end{equation*}
for $i = 1,2,3 \cdots$. This indicates that
\begin{equation*}
\int P(\rho) \,d \mu(\rho) = 0
\end{equation*}
for all $P(\rho) = \rho^j$ where $j = 0, 1, 2 \cdots$. Hence $\int P(\rho) \,d \mu(\rho) = 0$ for all polynomials. Moreover, the sequence $\{B^h(\Gamma_p^h)\}$ for $p = 1, 2, \cdots$ is absolute convergent. Therefore the compact support of this sequence is a bounded interval. Hence, by using Theorem~\ref{thm:brown} we have $\mu = 0$, which means $B^h(\Gamma_p^h) = 0$ for $p = 1, 2, \cdots, H[\Gamma]$. Similarly, we can obtain that $B^v(\Gamma_q^v) = 0$ for $q = 1, 2, \cdots, V[\Gamma]$. 
\end{proof}

We are now ready to provide the main result of this subsection, Theorem~\ref{thm:MUPIM}.
\begin{theorem}{\label{thm:MUPIM}}
If the measure induced by $\Gamma$ is the invariant measure of random walk $R$, then $|\Gamma_k^g| = \infty$ for $k = 1,2, \cdots, G[\Gamma]$ and $G[\Gamma] < \infty$, \ie $\Gamma$ must be the union of finitely many pairwise-coupled sets each with countably infinite cardinality.
\end{theorem}

\begin{proof}
Suppose that there exists a pairwise-coupled set $\Gamma_1^g \in \Gamma$ with $|\Gamma_1^g| < \infty$. From Lemma~\ref{lem:MUPIM} we have $B^h([\Gamma_1^g]_p^h) = 0$ for $p = 1, \dots, H[\Gamma_1^g]$ and $B^v([\Gamma_1^g]_q^v) = 0$ for $q = 1, \dots, V[\Gamma_1^g]$. Now, it follows from~\cite[Theorem~3]{chen2012invariant} that the measure induced by $\Gamma_1^g$ itself is an invariant measure of $R$. Since the measure induced by set $\Gamma$ is also an invariant measure, this contradicts the fact that the invariant measure is unique. Therefore, $|\Gamma_k^g| = \infty$ for $k = 1,2, \cdots, G[\Gamma]$.

Now, since $|\Gamma_k^g| = \infty$ and by Theorem~\ref{thm:singularity8}  there is at most one singularity on $Q$ in the unit square, it can be readily verified using Lemma~\ref{lem:MUPIM} that there exists a $(\rho, \sigma) \in \Gamma$ such that $B^h(\{(\rho, \sigma)\}) = 0$ or $B^v(\{(\rho, \sigma)\}) = 0$, which corresponds to the balance equations for states at the horizontal or at the vertical axis of the state space, respectively. In other words $(\rho,\sigma)\in H$ or $(\rho,\sigma)\in V$.

Finally, notice that there are at most a finite number of intersections between $Q$ and $H$ and between $Q$ and $V$, respectively. Therefore, we conclude that set $\Gamma$ must be the union of finitely many pairwise-coupled sets each with countably infinite cardinality.
\end{proof}

\subsection{Constraints on random walks}
In this section, we characterize the random walks for which the invariant measure may be an infinite sum of geometric terms. We will show that the existence of transitions to north, northeast or east plays an essential role in distinguishing such kind of random walks.
\begin{theorem}{\label{thm:propertyC}}
Let the invariant measure of random walk $R$ be $m(i,j) = \sum_{(\rho, \sigma) \in \Gamma} \alpha(\rho, \sigma) \rho^i \sigma^j$ with $|\Gamma| = \infty$. Then $p_{1,0} = p_{1,1} = p_{0,1} = 0$ and $\Gamma$ has a unique accumulation point in the origin.
\end{theorem}

\begin{proof}
We will demonstrate that in the absence of singularities it is not possible to have an invariant measure $m(i,j) = \sum_{(\rho, \sigma) \in \Gamma} \alpha(\rho, \sigma) \rho^i \sigma^j$ with $|\Gamma| = \infty$. More precisely, we show that $\Gamma$ must have a singularity of $Q$ as an accumulation point. The result of the theorem then follows, by Theorem~\ref{thm:singularity8}, that states the algebraic curve $Q$ has a singularity in $\bar{U}$ only if $p_{1,0} = p_{1,1} = p_{0,1} = 0$, in which case it is in the origin.

Suppose that $Q^+$ does not contain any singularities and that $m(i,j) = \sum_{(\rho, \sigma) \in \Gamma} \alpha(\rho, \sigma) \rho^i \sigma^j$ with $|\Gamma| = \infty$ is the invariant measure of $R$. In the remainder of this proof we will obtain a contradiction by showing that at least one of the terms $(\rho,\sigma)\in\Gamma$ will be outside the unit square and that the measure $m(i,j)$ can, therefore, not be finite. 

\begin{figure}
  \hfill
\subfigure[\label{fig:Qa}]
{
\begin{tikzpicture}[scale = 0.35]
\begin{axis}[%
view={0}{90},
width=4.52083333333333in,
height=4.52083333333333in,
scale only axis,
scale only axis,
xmin=0, xmax=1.4,
ymin=0, ymax=1.4,
xtick = {-0.5,0,0.5,1,1.4},
ytick = {-0.5,0.5,1,1.4}]

\addplot [draw=red, dashed, ultra thick] coordinates{ (0.202020202020202,0.353682214701604)(0.198384884176416,0.323232323232323)(0.197077895125914,0.282828282828283)
(0.202020202020202,0.246508855441227)(0.203210994526056,0.242424242424242)(0.242424242424242,0.203210994526056)(0.246508855441227,0.202020202020202)(0.282828282828283,0.197077895125914)
(0.323232323232323,0.198384884176416)(0.353682214701604,0.202020202020202)(0.363636363636364,0.203895429726532)(0.404040404040404,0.214494599816697)(0.444444444444444,0.225785196573688)
(0.484848484848485,0.236997817356168)(0.504404788147808,0.242424242424242)(0.525252525252525,0.250700917059662)(0.565656565656566,0.266180848525009)(0.606060606060606,0.28040464748359)
(0.613128657834267,0.282828282828283)(0.646464646464647,0.298405096645817)(0.686868686868687,0.315717377758959)(0.705478529499738,0.323232323232323)(0.727272727272727,0.334879782200734)
(0.767676767676768,0.355000115222371)(0.786170830398225,0.363636363636364)(0.808080808080808,0.376913130462884)(0.848484848484849,0.399657695376722)(0.856667825816326,0.404040404040404)
(0.888888888888889,0.426152914840109)(0.917591296514448,0.444444444444444)(0.929292929292929,0.453964973403925)(0.969690941847002,0.484848484848485)(0.96969696969697,0.484854384740536)
(1.01010101010101,0.522600209342491)(1.01309284479204,0.525252525252525)(1.04839442169988,0.565656565656566)(1.05050505050505,0.568865091941115)(1.07584234204985,0.606060606060606)
(1.09090909090909,0.636320777503951)(1.09610259029538,0.646464646464647)(1.10924344886468,0.686868686868687)(1.11585119923513,0.727272727272727)(1.11615563320778,0.767676767676768)
(1.11040146355617,0.808080808080808)(1.09884567781367,0.848484848484849)(1.09090909090909,0.867567246365329)(1.08131783270928,0.888888888888889)(1.05778644774189,0.929292929292929)
(1.05050505050505,0.939687761385969)(1.02741480688696,0.96969696969697)(1.01010101010101,0.989122644265788)(0.989122644265788,1.01010101010101)(0.96969696969697,1.02741480688696)
(0.939687761385969,1.05050505050505)(0.929292929292929,1.05778644774189)(0.888888888888889,1.08131783270928)(0.867567246365329,1.09090909090909)(0.848484848484849,1.09884567781367)
(0.808080808080808,1.11040146355617)(0.767676767676768,1.11615563320778)(0.727272727272727,1.11585119923513)(0.686868686868687,1.10924344886468)(0.646464646464647,1.09610259029538)
(0.636320777503952,1.09090909090909)(0.606060606060606,1.07584234204985)(0.568865091941115,1.05050505050505)(0.565656565656566,1.04839442169988)(0.525252525252525,1.01309284479204)
(0.522600209342491,1.01010101010101)(0.484854384740536,0.96969696969697)(0.484848484848485,0.969690941847002)(0.453964973403925,0.929292929292929)(0.444444444444444,0.917591296514448)
(0.426152914840109,0.888888888888889)(0.404040404040404,0.856667825816326)(0.399657695376722,0.848484848484848)(0.376913130462884,0.808080808080808)(0.363636363636364,0.786170830398225)
(0.355000115222371,0.767676767676768)(0.334879782200734,0.727272727272727)(0.323232323232323,0.705478529499738)(0.315717377758959,0.686868686868687)(0.298405096645817,0.646464646464647)
(0.282828282828283,0.613128657834268)(0.28040464748359,0.606060606060606)(0.266180848525009,0.565656565656566)(0.250700917059662,0.525252525252525)(0.242424242424242,0.504404788147808)
(0.236997817356168,0.484848484848485)(0.225785196573688,0.444444444444444)(0.214494599816697,0.404040404040404)(0.203895429726532,0.363636363636364)(0.202020202020202,0.353682214701604)
(NaN,NaN)};
\addplot [draw=red, ultra thick] coordinates{(0.282828282828283,0.197077895125914)
(0.323232323232323,0.198384884176416)(0.353682214701604,0.202020202020202)(0.363636363636364,0.203895429726532)(0.404040404040404,0.214494599816697)(0.444444444444444,0.225785196573688)
(0.484848484848485,0.236997817356168)(0.504404788147808,0.242424242424242)(0.525252525252525,0.250700917059662)(0.565656565656566,0.266180848525009)(0.606060606060606,0.28040464748359)
(0.613128657834267,0.282828282828283)(0.646464646464647,0.298405096645817)(0.686868686868687,0.315717377758959)(0.705478529499738,0.323232323232323)(0.727272727272727,0.334879782200734)
(0.767676767676768,0.355000115222371)};
\addplot [draw=red, ultra thick] coordinates{(0.767676767676768,1.11615563320778)(0.727272727272727,1.11585119923513)(0.686868686868687,1.10924344886468)(0.646464646464647,1.09610259029538)
(0.636320777503952,1.09090909090909)(0.606060606060606,1.07584234204985)(0.568865091941115,1.05050505050505)(0.565656565656566,1.04839442169988)(0.525252525252525,1.01309284479204)
(0.522600209342491,1.01010101010101)(0.484854384740536,0.96969696969697)(0.484848484848485,0.969690941847002)(0.453964973403925,0.929292929292929)(0.444444444444444,0.917591296514448)
(0.426152914840109,0.888888888888889)(0.404040404040404,0.856667825816326)(0.399657695376722,0.848484848484848)(0.376913130462884,0.808080808080808)(0.363636363636364,0.786170830398225)
(0.355000115222371,0.767676767676768)(0.334879782200734,0.727272727272727)(0.323232323232323,0.705478529499738)(0.315717377758959,0.686868686868687)(0.298405096645817,0.646464646464647)
};
\addplot [
color=black,
solid
]
coordinates{
 (0,1)(1,1) 
};

\addplot [
color=black,
solid
]
coordinates{
 (1,0)(1,1) 
};

\addplot[color=black,dashed,thick,mark=*, mark options={fill=white}] 
    coordinates {
         (0.2,0.3)
         (0.2,0)
        };

\addplot[color=black,dashed,thick,mark=*, mark options={fill=white}] 
    coordinates {
         (0.3, 0.2)
         (0, 0.2)
        };

\addplot[color=black,dashed,thick,mark=*, mark options={fill=white}] 
    coordinates {
         (0.76, 1.11)
         (0, 1.11)
        };

\addplot[color=black,dashed,thick,mark=*, mark options={fill=white}] 
    coordinates {
         (1.11, 0.76)
         (1.11, 0)
        };
\addplot[color=black,solid,thick,mark=*, mark options={fill=white}] 
    coordinates {
         (0.76, 1.11)
         (0.76, 0.35)
        };
\addplot[color=black,solid,thick,mark=*, mark options={fill=white}] 
    coordinates {
         (0.3,0.2)
        (0.3,0.64)
        };       
\node [above right] at (axis cs:  0,  1.11) {\Huge $\scriptstyle y_t$};
\node [above right] at (axis cs:  1.11, 0) {\Huge $\scriptstyle x_r$};
\node [below right] at (axis cs:  0, 0.2) {\Huge $\scriptstyle y_b$};
\node [above right] at (axis cs:  0.2, 0) {\Huge $\scriptstyle x_l$};
\node at (axis cs:  0.25,  0.25) {\Huge $\scriptstyle Q_{00}$};
\node at (axis cs:  0.7,  0.33) {\Huge $\scriptstyle Q_{10}$};
\node at (axis cs:  1,  1) {\Huge $\scriptstyle Q_{11}$};
\node at (axis cs:  0.4,  0.77) {\Huge $\scriptstyle Q_{01}$};
\node at (axis cs:  0.25,  0.4) {\Huge $Q_l$};
\node at (axis cs:  0.5,  0.5) {\Huge $Q_c$};
\node at (axis cs:  0.85,  0.6) {\Huge $Q_r$};
\end{axis}
\end{tikzpicture}
} 
\subfigure[\label{fig:Qb}]
{
\begin{tikzpicture}[scale = 0.35]
\begin{axis}[%
view={0}{90},
width=4.52083333333333in,
height=4.52083333333333in,
scale only axis,
scale only axis,
xmin=0, xmax=1.4,
ymin=0, ymax=1.4,
xtick = {-0.5,0,0.5,1,1.4},
ytick = {-0.5,0.5,1,1.4}]

\addplot [draw=red, ultra thick] coordinates{ (0.202020202020202,0.353682214701604)(0.198384884176416,0.323232323232323)(0.197077895125914,0.282828282828283)(0.202020202020202,0.246508855441227)(0.203210994526056,0.242424242424242)(0.242424242424242,0.203210994526056)(0.246508855441227,0.202020202020202)(0.282828282828283,0.197077895125914)(0.323232323232323,0.198384884176416)(0.353682214701604,0.202020202020202)(0.363636363636364,0.203895429726532)(0.404040404040404,0.214494599816697)(0.444444444444444,0.225785196573688)(0.484848484848485,0.236997817356168)(0.504404788147808,0.242424242424242)(0.525252525252525,0.250700917059662)(0.565656565656566,0.266180848525009)(0.606060606060606,0.28040464748359)(0.613128657834267,0.282828282828283)(0.646464646464647,0.298405096645817)(0.686868686868687,0.315717377758959)(0.705478529499738,0.323232323232323)(0.727272727272727,0.334879782200734)(0.767676767676768,0.355000115222371)(0.786170830398225,0.363636363636364)(0.808080808080808,0.376913130462884)(0.848484848484849,0.399657695376722)(0.856667825816326,0.404040404040404)(0.888888888888889,0.426152914840109)(0.917591296514448,0.444444444444444)(0.929292929292929,0.453964973403925)(0.969690941847002,0.484848484848485)(0.96969696969697,0.484854384740536)(1.01010101010101,0.522600209342491)(1.01309284479204,0.525252525252525)(1.04839442169988,0.565656565656566)(1.05050505050505,0.568865091941115)(1.07584234204985,0.606060606060606)(1.09090909090909,0.636320777503951)(1.09610259029538,0.646464646464647)(1.10924344886468,0.686868686868687)(1.11585119923513,0.727272727272727)(1.11615563320778,0.767676767676768)(1.11040146355617,0.808080808080808)(1.09884567781367,0.848484848484849)(1.09090909090909,0.867567246365329)(1.08131783270928,0.888888888888889)(1.05778644774189,0.929292929292929)(1.05050505050505,0.939687761385969)(1.02741480688696,0.96969696969697)(1.01010101010101,0.989122644265788)(0.989122644265788,1.01010101010101)(0.96969696969697,1.02741480688696)(0.939687761385969,1.05050505050505)(0.929292929292929,1.05778644774189)(0.888888888888889,1.08131783270928)(0.867567246365329,1.09090909090909)(0.848484848484849,1.09884567781367)(0.808080808080808,1.11040146355617)(0.767676767676768,1.11615563320778)(0.727272727272727,1.11585119923513)(0.686868686868687,1.10924344886468)(0.646464646464647,1.09610259029538)(0.636320777503952,1.09090909090909)(0.606060606060606,1.07584234204985)(0.568865091941115,1.05050505050505)(0.565656565656566,1.04839442169988)(0.525252525252525,1.01309284479204)(0.522600209342491,1.01010101010101)(0.484854384740536,0.96969696969697)(0.484848484848485,0.969690941847002)(0.453964973403925,0.929292929292929)(0.444444444444444,0.917591296514448)(0.426152914840109,0.888888888888889)(0.404040404040404,0.856667825816326)(0.399657695376722,0.848484848484848)(0.376913130462884,0.808080808080808)(0.363636363636364,0.786170830398225)(0.355000115222371,0.767676767676768)(0.334879782200734,0.727272727272727)(0.323232323232323,0.705478529499738)(0.315717377758959,0.686868686868687)(0.298405096645817,0.646464646464647)(0.282828282828283,0.613128657834268)(0.28040464748359,0.606060606060606)(0.266180848525009,0.565656565656566)(0.250700917059662,0.525252525252525)(0.242424242424242,0.504404788147808)(0.236997817356168,0.484848484848485)(0.225785196573688,0.444444444444444)(0.214494599816697,0.404040404040404)(0.203895429726532,0.363636363636364)(0.202020202020202,0.353682214701604)(NaN,NaN)};
\addplot [
color=black,
solid
]
coordinates{
 (0,1)(1,1) 
};

\addplot [
color=black,
solid
]
coordinates{
 (1,0)(1,1) 
};

\addplot[color=black,solid,thick,mark=*, mark options={fill=white}] 
    coordinates {
         (0.95, 0.47)
         (0.235, 0.47)
        };
\addplot[color=black,solid,thick,mark=*, mark options={fill=white}] 
    coordinates {
         (0.235, 0.47)
         (0.235, 0.215)
        };
\addplot[color=black,dashed,thick,mark=*, mark options={fill=white}] 
    coordinates {
         (0.235, 0.215)
         (0.4, 0.215)
        };
\addplot[color=black,dashed,thick,mark=*, mark options={fill=white}] 
    coordinates {
         (0.4, 0.215)
         (0.4, 0.84)
        };
\addplot[color=black,dashed,thick,mark=*, mark options={fill=white}] 
    coordinates {
         (0.4, 0.84)
         (1.1, 0.84)
        };
\addplot[color=black,solid,thick,mark=*, mark options={fill=white}] 
    coordinates {
         (1.1, 0.84)
         (1.1, 0.665)
        };
\addplot[color=black,solid,thick,mark=*, mark options={fill=white}] 
    coordinates {
         (1.1, 0.665)
         (0.31, 0.665)
        };
\addplot[color=black,solid,thick,mark=*, mark options={fill=white}] 
    coordinates {
         (0.31, 0.665)
         (0.31, 0.2)
        };
\addplot[color=black,solid,thick,mark=*, mark options={fill=white}] 
    coordinates {
         (0.235, 0)
         (0.235, 0)
        };
\addplot[color=black,solid,thick,mark=*, mark options={fill=white}] 
    coordinates {
         (1.1, 0)
         (1.1, 0)
        };
\node [above] at (axis cs:  0.235, 0) {\Huge $\scriptstyle \rho_{1,L(1)}$};
\node [above] at (axis cs:  1.1, 0) {\Huge $\scriptstyle \rho_{2,1}$};
\node at (axis cs:  0.5,  0.47) {\Huge $\scriptstyle \Gamma_1$};
\node at (axis cs:  0.7,  0.84) {\Huge $\scriptstyle \Gamma_2$};
\node at (axis cs:  0.6,  0.665) {\Huge $\scriptstyle \Gamma_3$};
\end{axis}
\end{tikzpicture}
} 
  \hfill{}
   \caption{~\subref{fig:Qa} Different partition of $Q^{+}$ for the random walk in Figure~\ref{fig:2c}.~\subref{fig:Qb} Example of $\Gamma_1, \Gamma_2$ and $\Gamma_3$. {\label{fig:middleQ}}}
\end{figure}
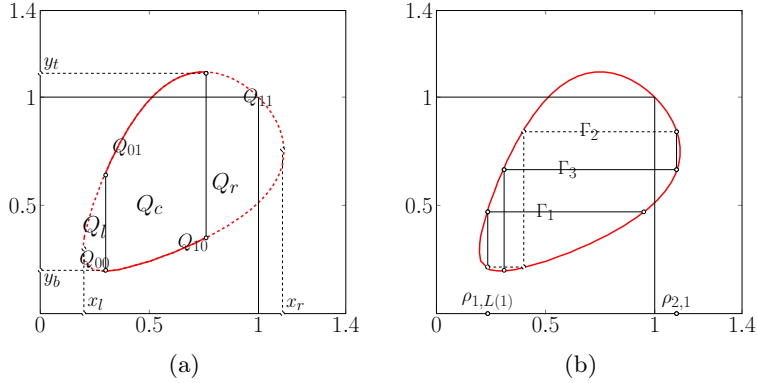

Consider $E(x,y) = Q(e^x, e^y) = 0$ for $(x,y) \in \mathbb{R}^2$. The set of pairs $(x,y)$ satisfying $E(x,y) = 0$ is denoted by $E$. It has been shown that $E(x,y) < 0$ is convex in~\cite{latouche2014product, miyazawa2009tail}, examples of $E(x,y) = 0$ can be found in~\cite{latouche2014product, miyazawa2009tail}. It is readily verified that logarithm serves as a bijection between $(x,y)$ satisfying $E(x,y) = 0$ and geometric term $(\rho,\sigma)$ satisfying $Q(\rho,\sigma) = 0$. Moreover, curve $E$ is bounded by $[\log x_l, \log x_r] \times [\log y_b, \log y_t]$, \ie in terms of the branching points of $Q$. Since there is no singularity on $Q$, we have $0 < x_l < 1$ and $0 < y_b < 1$. Hence, we have $\log x_l > -\infty$ and $\log y_b > -\infty$.

If $\Gamma = \{(\rho_1, \sigma_1), (\rho_2, \sigma_2), (\rho_3, \sigma_3), \cdots\}$ is a pairwise-coupled set on $Q$, then $\tilde{\Gamma} = \{(\log{\rho_1}, \log{\sigma_1})$, $(\log{\rho_2}, \log{\sigma_2})$, $(\log{\rho_3}, \log{\sigma_3}), \cdots \}$ is pairwise-coupled on $E$. We have $\rho < 1$ and $\sigma < 1$ for any $(\rho, \sigma) \in \Gamma$. Therefore, we have $x < 0$ and $y < 0$ for any $(x,y) \in \tilde{\Gamma}$.

We now construct a partition $\{\Gamma_1,\dots,\Gamma_K\}$ of $\Gamma$, where the elements of $\Gamma_i$ are denoted as $\Gamma_i=\{(\rho_{i,1},\sigma_{i,1}),\dots,(\rho_{i,L(i)},\sigma_{i,L(i)})\}$ and each $\Gamma_i$ satisfies
 \begin{equation} \label{eq:thmpropertyCstair}
 \begin{IEEEeqnarraybox}[][c]{rClrCl}
 \rho_{i,1} &>& \rho_{i,2},\quad & \sigma_{i,1} &=& \sigma_{i,2}, \\
 \rho_{i,2} &=& \rho_{i,3},\quad & \sigma_{i,2} &>& \sigma_{i,3}, \\ 
 \rho_{i,3} &>& \rho_{i,4},\quad & \sigma_{i,3} &=& \sigma_{i,4}, \\
 &\vdots& & &\vdots& \\
 \rho_{i,L(i)-1} &>& \rho_{i,L(i)},\quad & \sigma_{i,L(i)-1} &=& \sigma_{i,L(i)},
 \end{IEEEeqnarraybox}
 \end{equation}

In addition the partition $\{\Gamma_1,\dots,\Gamma_K\}$ is maximal in the sense that no $\Gamma_i\cup\Gamma_j$, $i\neq j$ satisfies~\eqref{eq:thmpropertyCstair}. Let $\{\tilde{\Gamma}_1, \tilde{\Gamma}_2, \tilde{\Gamma_3}, \cdots\}$ denote the corresponding partition of $\tilde{\Gamma}$. The partition is illustrated in Figure~\ref{fig:Qb}.

Consider set $\tilde{\Gamma}$ on $E$. It can be readily verified that $|\tilde{\Gamma}_k| < \infty$ where $k = 1,2, \cdots, K$ because $E(x,y) < 0$ where $x < 0$ and $y < 0$ is a convex set and $\log x_l > -\infty$ and $\log y_b > -\infty$. Therefore, we conclude that $|\Gamma_k| < \infty$ where $k = 1,2, \cdots, K$ as well.

Assume $x_b \leq x_t$, the case where $x_b > x_b$ can be analyzed similarly. To simplify the presentation we introduce additional notation. Let $\{Q_l, Q_c, Q_r\}$ denote a partition of $Q$, where
\begin{equation*}
\begin{aligned}
Q_l &= \left\{ (x,y)\in Q\; \middle|\; x \leq x_b \right\}, \\
Q_c &= \left\{ (x,y)\in Q\; \middle|\; x_b <x \leq x_t \right\}, \\
Q_r &= \left\{ (x,y)\in Q\; \middle|\; x > x_t \right\}.
\end{aligned}
\end{equation*}

Next, we prove $K<\infty$. More precisely, we will show that $K\leq 2$. Without loss of generality, we assume $K=3$ and $|\Gamma_i| \geq 2$. Observe that $\{\Gamma_1,\Gamma_2,\Gamma_3\}$ forms a pairwise-coupled set. Therefore, we must have $\rho_{1,L(1)} = \rho_{2,L(2)}$  and $\rho_{2, 1} = \rho_{3,1}$. This is illustrated in Figure~\ref{fig:Qb}. From the monotonicity of $Q_{10}$ and the structure of the partition $\{\Gamma_1,\dots,\Gamma_K\}$ it follows that $\rho_{1,L(1)} \in Q_l$.  Similarly it follows that $\rho_{2,1} \in Q_r$ and $\rho_{3,1} \in Q_r$.  Moreover, one of $\rho_{2,1}$ and $\rho_{3,1}$ must be on $Q_{11}$. Note that from Lemma~\ref{lem:2389} we have $y_t \geq 1$ and $x_r \geq 1$. Together with the monotonicity of $Q_{11}$ from Lemma~\ref{lem:monotonicity} this leads to the conclusion that $Q_{11}$ is outside of $U$. This implies that one of the geometric measures in $\Gamma$ is outside $U$ and to a contradiction to the fact the induced measure is finite. 

We have shown that under the assumption that $Q^+$ does not contain singularities, $\Gamma$ can be partitioned into a finite number of sets each with a finite number of elements.
\end{proof}
%

\subsection{Constraints on the coefficients}
The last section was devoted to finding the constraint to the random walk in which the pairwise-coupled set with infinite cardinality could be obtained. In this section, we show that it is necessary to have a geometric term with negative coefficient in the linear combination of infinite geometric terms.

\begin{theorem}\label{thm:negative}
If the invariant measure of random walk $R$ is $m(i,j) = \sum_{(\rho, \sigma) \in \Gamma} \alpha(\rho, \sigma)\rho^i\sigma^j$, where $\Gamma\subset Q^{+}_{U}$, $|\Gamma| = \infty$ and $\alpha(\rho,\sigma)\in \mathbb{R} \backslash \{0\}$, then at least one $\alpha(\rho, \sigma)$ is negative.
\end{theorem}
\begin{proof}
We know from Theorem~\ref{thm:propertyC} that $p_{1,0} = p_{1,1} = p_{0,1} = 0$ in random walk $R$. Notice that there is at least one pairwise-coupled set from $\Gamma$ that contains countably infinitely many geometric terms due to Theorem~\ref{thm:MUPIM}. Without loss of generality, we consider set $\Gamma$ which only contains a single pairwise-coupled set which is of the form given in~\eqref{eq:pairwisexample} and assume that $\alpha(\rho_1,\sigma_1)=1$. Since measure induced by $\Gamma$ is the invariant measure, it follows from Theorem~\ref{thm:MUPIM} that 
$B^h(\Gamma_p^h) = 0$ and $B^v(\Gamma_q^v) = 0$ for all $p,q \in \{1,2,3, \cdots\}$. Note that $B^h\{(\rho_1, \sigma_1), (\rho_2, \sigma_2)\} = 0$ indicates that $\alpha(\rho_2, \sigma_2)$ is uniquely determined by
\begin{equation*}
\alpha(\rho_2, \sigma_2) = - \frac{T_1}{T_2} \alpha(\rho_1,\sigma_1),
\end{equation*}
where
\begin{equation*}
T_i = (1 - \frac{1}{\rho_i}) h_1 + (1 - \rho_i)h_{-1} + \sum_{s = -1}^{1} p_{s,1} - \sigma_i(\sum_{s = -1}^{1} \rho_i^{-s}p_{s,-1}).
\end{equation*}
Next, $\alpha(\rho_3, \sigma_3)$ follows from $B^v\{(\rho_2, \sigma_2), (\rho_3, \sigma_3)\} = 0$. In similar fashion, for $k \in \{1,2,3, \cdots\}$, 
\begin{equation*}
\alpha(\rho_{2k}, \sigma_{2k}) = - \frac{T_{2k - 1}}{T_{2k}}  \alpha(\rho_{2k-1},\sigma_{2k-1}) \quad \text{where} \quad \rho_{2k} = \rho_{2k - 1}, \sigma_{2k} < \sigma_{2k - 1}.
\end{equation*}
The following two facts allow us to show that there exists a positive integer $N$ such that $\frac{T_{2k - 1}}{T_{2k}} > 0$ when $k > N$. First, we know $p_{1,0} = p_{1,1} = p_{0,1} = 0$ and $\lim_{k \rightarrow \infty} \rho_{k} = 0$ from Theorem~\ref{thm:propertyC}. Secondly, the ergodic random walk with no drift to northeast requires $h_1  +  p_{1,-1} \neq 0$. Note
\begin{equation*}
\frac{T_{2k - 1}}{T_{2k}} = \frac{(1 - \frac{1}{\rho_{2k-1}}) h_1 + (1 - \rho_{2k - 1}) h_{-1} + p_{-1,1} - \sigma_{2k - 1}(\sum_{s = -1}^{1} \rho_{2k-1}^{-s}p_{s,-1})}{(1 - \frac{1}{\rho_{2k}}) h_1 + (1 - \rho_{2k}) h_{-1} + p_{-1,1} - \sigma_{2k}(\sum_{s = -1}^{1} \rho_{2k}^{-s}p_{s,-1})},
\end{equation*}
By using L'Hospital's rule, we can conclude
\begin{equation*}
\lim_{k \rightarrow \infty} \frac{T_{2k - 1}}{T_{2k}} = \begin{cases} \frac{h_{1} + \sigma_{2k - 1} p_{1,-1}}{h_{1} + \sigma_{2k} p_{1,-1}}, &\mbox{if } h_1 \neq 0, \quad p_{1,-1} \neq 0; \\ 
1, & \mbox{if } h_1 \neq 0, \quad p_{1,-1} = 0; \\
\frac{\sigma_{2k - 1}}{\sigma_{2k}}, & \mbox{if } h_1 = 0, \quad p_{1,-1} \neq 0. \end{cases} 
\end{equation*}
The non-negativity of $\frac{T_{2k - 1}}{T_{2k}}$ when $k$ is large enough completes the proof.
\end{proof}

\section{Example: $2 \times 2$ Switch}{\label{sec:switch}}
In this section we provide an example of a random walk for which the invariant measure satisfies all conditions that have been obtained in the current paper. In particular, we consider the $2 \times 2$ switch which has been studied by Boxma and Van Houtum in~\cite{boxma1993compensation}.

\begin{figure}
\begin{center}
\begin{tikzpicture}[>=latex]
\draw (0,0) -- ++(2cm,0) -- ++(0,-1.5cm) -- ++(-2cm,0);
\foreach \i in {1,...,4}
  \draw (2cm-\i*10pt,0) -- +(0,-1.5cm);

\draw (2.75,-0.75cm) circle [radius=0.75cm];

\draw[->] (3.5,-0.75) -- +(30pt,0);
\draw (0,-0.75) -- +(-180pt,0) node[left] {$r_2$};
\node at (2.75,-0.75cm) {$2$};

\draw (0,3) -- ++(2cm,0) -- ++(0,-1.5cm) -- ++(-2cm,0);
\foreach \i in {1,...,4}
  \draw (2cm-\i*10pt,3) -- +(0,-1.5cm);

\draw (2.75,2.25cm) circle [radius=0.75cm];

\draw[->] (3.5,2.25) -- +(30pt,0);
\draw (0,2.25) -- +(-180pt,0) node[left] {$r_1$};
\node at (2.75,2.25cm) {$1$};
\draw (-4,2.25) -- (0,-0.75);
\draw (-4,-0.75) -- (0,2.25);
\draw[->] (-6,2.25) --(-5,2.25);
\draw[->] (-6,-0.75) --(-5,-0.75);
\draw[->] (-4,2.25) --(-2,2.25);
\draw[->] (-4,-0.75) --(-2,-0.75);
\draw[->] (-4,2.25) --(-3,1.5);
\draw[->] (-4,-0.75) --(-3,0);

\node at (-2,2.5cm) {$t_{11}$};
\node at (-3.5,1.5cm) {$t_{12}$};
\node at (-2,-1cm) {$t_{22}$};
\node at (-3.5,0cm) {$t_{21}$};
\end{tikzpicture}
\caption{The $2 \times 2$ switch. \label{fig:switch}}
\end{center}
\end{figure}
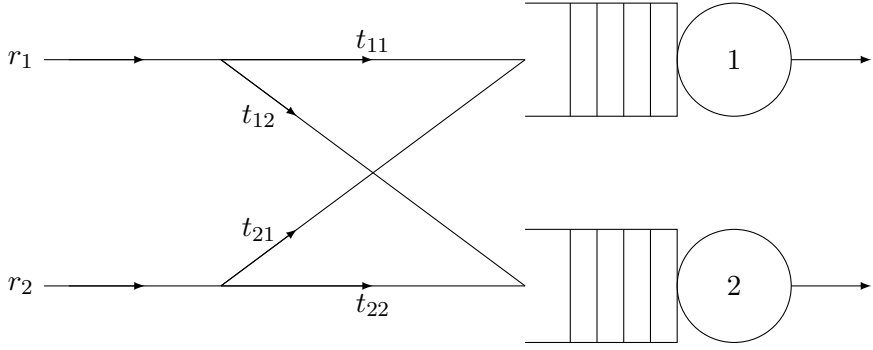

A $2 \times 2$ buffer switch has two input and two output ports. Such a switch is modeled as a discrete time queueing system with two parallel servers and two type of arriving jobs(see Figure~\ref{fig:switch}). Jobs of type $i$, $i = 1,2$, are assumed to arrive according to a Bernoulli stream with rate $r_i$, $0 < r_i \leq 1$. This means that at every time unit the number of arriving jobs of type $i$ is one with probability $r_i$ and zero with probability $1 - r_i$. Jobs always arrive at the beginning of a time unit, and once a job of type $i$ has arrived, it joins the queue at server $j$ with probability $t_{ij}$, $t_{ij} > 0$ for $j = 1,2$, and $t_{i,1} + t_{i,2} = 1$. Jobs that have arrive at the beginning of a time unit are immediately candidates for service. A server serves exactly one job per time unit, if one is present. We assume the system is stable. 

\begin{figure}
\hfill
\subfigure[
{\label{fig:2by2A}}
]
{\hfill
\begin{tikzpicture}[scale=0.75]
\tikzstyle{axes}=[very thin] \tikzstyle{trans}=[very thick,->]
   \draw[axes] (0,0)  -- node[at end, below] {$\scriptstyle \rightarrow i$} (6,0); 
   \draw[axes] (0,0) -- node[at end, left] {$\scriptstyle {\uparrow} {j}$} (0,5);
   \draw[trans] (0,0) to node[below] {$\scriptstyle p_{1,-1}$} (1,0);
   \draw[trans] (0,0) to node[left] {$\scriptstyle p_{-1,1}$} (0,1);
   \draw[trans] (4,0) to node[at end, anchor = north] {$\scriptstyle p_{-1,0} + p_{-1,-1}$} (3,0);
   \draw[trans] (4,0) to node[below right] {$\scriptstyle p_{1,-1}$} (5,0);
   \draw[trans] (4,0) to node[at end, anchor = south east] {$\scriptstyle p_{-1,1}$} (3,1);
   \draw[trans] (0,3) to node[left] {$\scriptstyle p_{0,-1}+p_{-1,-1}$} (0,2);
   \draw[trans] (0,3) to node[left] {$\scriptstyle p_{-1,1}$} (0,4);
   \draw[trans] (0,3) to node[at end, anchor = north] {$\scriptstyle p_{1,-1}$} (1,2);
   \draw[trans] (4,3) to node[at end, anchor = south east] {$\scriptstyle p_{-1,1}$} (3,4);
   \draw[trans] (4,3) to node[at end, anchor = east] {$\scriptstyle p_{-1,0}$} (3,3);
   \draw[trans] (4,3) to node[at end, anchor = north east] {$\scriptstyle p_{-1,-1}$} (3,2);
   \draw[trans] (4,3) to node[at end, anchor = north] {$\scriptstyle p_{0,-1}$} (4,2);
   \draw[trans] (4,3) to node[at end, anchor = north west] {$\scriptstyle p_{1,-1}$} (5,2);  
\end{tikzpicture}
\hfill{}
} 
\subfigure[
 {\label{fig:2by2B}}
 ] 
 {
 \input{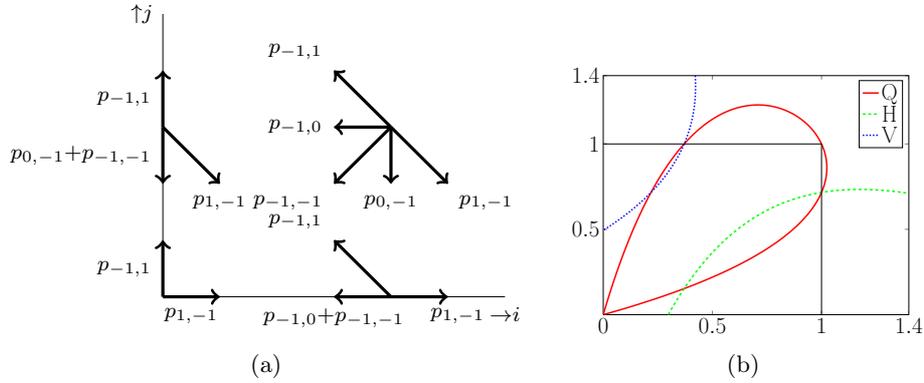}
 }
\hfill{}
\caption{The $2 \times 2$ switch.~\subref{fig:2by2A} Transition diagram.~\subref{fig:2by2B} Balance equations for the $2\times 2$ switch with $r_1 = 0.8$, $r_2 = 0.9$, $t_{11} = 0.3$, $t_{12} = 0.7$, $t_{21} = 0.6$, $t_{22} = 0.4$.\label{fig:2by2}}
\end{figure}

We now describe the $2 \times 2$ switch by a random walk in the quarter-plane with states $(i,j)$, where $i$ and $j$ denote the numbers of waiting jobs at server $1$ and server $2$, respectively, at the beginning of a time unit. For a state $(i,j)$ in the interior of the state space, we only have transitions to the neighboring state $(i+s, j+t)$ with $s,t \in \{-1,0,1\}$ and $s + t \leq 0$. The corresponding transition probabilities $p_{s,t}$ are equal to 
\begin{align*}
p_{1,-1} &= r_1 r_2 t_{11} t_{12}, \\
p_{0,0} &= r_1 r_2(t_{11}t_{22} + t_{12} t_{21}), \\
p_{-1,1} &= r_1 r_2 t_{12} t_{22}, \\
p_{0,-1} &= r_1 (1 - r_2)t_{11} + r_2(1 - r_1)t_{21}, \\
p_{-1,0} &= r_1(1 - r_2)t_{12} + r_2(1 - r_1) t_{22}, \\
p_{-1,-1} &= (1 - r_1)(1- r_2).
\end{align*}
Each transition probability for the states at the boundaries can be written as a sum of the probabilities $p_{s,t}$. In Figure~\ref{fig:2by2A} all transition probabilities, except those for the transition from a state to itself, are illustrated. In Figure~\ref{fig:2by2B} the algebraic curves for $Q$, $H$ and $V$ are depicted.

It has been shown in~\cite{boxma1993compensation} that the invariant measure for the $2 \times 2$ switch is the sum of two alternating series of geometric terms, both of which have infinite cardinality and are pairwise-coupled.

\section{Conclusion}\label{sec:conclusion}
In this paper, we have obtained necessary conditions for the invariant measure of a random walk to be an infinite sum of geometric terms. These conditions hold under a mild regularity condition on the structure of these terms. Firstly, Theorem~\ref{thm:productterm} says that each geometric term in the linear combination must satisfy the interior balance equation. Secondly, Theorem~\ref{thm:MUPIM} indicates that only a measure induced by finitely many pairwise-coupled sets each with countably infinite cardinality may yield an invariant measure. Thirdly, Theorem~\ref{thm:propertyC} shows that the invariant measure may be an infinite sum of geometric terms only if there are no transitions to the North, Northeast or East. Finally, Theorem~\ref{thm:negative} requires at least one of the coefficients in the linear combination to be negative. 

A measure $m$ that satisfied the above necessary conditions may indeed be the invariant measure of a random walk. For example, the invariant measure for the $2 \times 2$ switch problem satisfies all the above conditions~\cite{boxma1993compensation}.


\section{Acknowledgment}
The authors thank Anton A. Stoorvogel for useful discussions. The authors also wish to thank an anonymous reviewer for useful comments that have helped to improve the presentation of the paper.  Yanting Chen acknowledges support through a CSC scholarship [No.2008613008]. This work is partly supported by the Netherlands Organization for Scientific Research (NWO) grant $612.001.107$.

\bibliographystyle{plain}
\bibliography{bibfile}
\end{document}